\documentclass{amsart}
\usepackage[hyphens]{url}
\usepackage{mathtools, microtype, mathrsfs, xfrac, hyperref, amssymb, enumerate, xspace, stmaryrd}
\usepackage{bbm}
\usepackage[utf8]{inputenc}
\usepackage{tikz-cd}
\usepackage{pdflscape}
\usepackage{bm}

\usepackage[backend=biber, style=alphabetic]{biblatex}
\usepackage{todonotes}

\numberwithin{equation}{section}

\numberwithin{subsection}{section}

\allowdisplaybreaks[1]

\AtBeginDocument{%
   \def\MR#1{}
}

\newtheorem*{namedtheorem}{\theoremname}
\newcommand{\theoremname}{testing}

\newtheorem{theorem}{Theorem}[section]
\newtheorem{proposition}[theorem]{Proposition}
\newtheorem{corollary}[theorem]{Corollary}
\newtheorem{lemma}[theorem]{Lemma}
\newtheorem*{theorem*}{Theorem}

\theoremstyle{definition}
\newtheorem{definition}[theorem]{Definition}
\newtheorem{notation}[theorem]{Notation}
\newtheorem{convention}[theorem]{Convention}
\newtheorem{conjecture}[theorem]{Conjecture}
\newtheorem{example}[theorem]{Example}

\newtheorem{remark}[theorem]{Remark}

\newtheorem*{question*}{Question}

\theoremstyle{remark}

\renewcommand{\mathcal}{\mathscr}

\newcommand\cA{\mathcal{A}} \newcommand\cB{\mathcal{B}}
\newcommand\cC{\mathcal{C}} 
 
\newcommand\cG{\mathcal{G}} \newcommand\cH{\mathcal{H}}

\newcommand\cO{\mathcal{O}} \newcommand\cP{\mathcal{P}}
 
\newcommand\cS{\mathcal{S}} \newcommand\cT{\mathcal{T}}
 
\newcommand\cW{\mathcal{W}}

\newcommand\CC{\mathbb{C}} 
 
\newcommand\GG{\mathbb{G}}

 \newcommand\NN{\mathbb{N}}
 \newcommand\PP{\mathbb{P}}
\newcommand\QQ{\mathbb{Q}} 
\renewcommand\SS{\mathbb{S}}

 \newcommand\ZZ{\mathbb{Z}}

 \newcommand\rF{\mathrm{F}}
 
\newcommand\rI{\mathrm{I}} 
 
 \newcommand\rN{\mathrm{N}}
 \newcommand\rP{\mathrm{P}}

\newcommand\rmq{\mathrm{q}}

 \newcommand\frm{\mathfrak{m}}
 
 \newcommand\fro{\mathfrak{o}}

\newcommand\arr{\ifinner\to\else\longrightarrow\fi}

\newcommand\arrto{\ifinner\mapsto\else\longmapsto\fi}

\renewcommand\H{\operatorname{H}}

\newcommand{\eqdef}{\mathrel{\smash{\overset{\mathrm{\scriptscriptstyle def}} =}}}

\def\displaytimes_#1{\mathrel{\mathop{\times}\limits_{#1}}}

\def\displayotimes_#1{\mathrel{\mathop{\bigotimes}\limits_{#1}}}

\renewcommand\hom{\operatorname{Hom}}

\newcommand\ext{\operatorname{Ext}}

\newcommand\aut{\operatorname{Aut}}

\newcommand\tor{\operatorname{Tor}}

\newcommand\pic{\operatorname{Pic}}

\newcommand\spec{\operatorname{Spec}}

\newcommand\rk{\operatorname{rk}}

\newcommand\id{\mathrm{id}}

\newcommand\indlim{\varinjlim}

\renewcommand\projlim{\varprojlim}

\newcommand{\underaut}{\mathop{\underline{\mathrm{Aut}}}\nolimits}

\newlength{\ignora}

\renewcommand{\setminus}{\smallsetminus}

\renewcommand\projlim{\varprojlim}

\newcommand{\GL}{\mathrm{GL}}

\newcommand{\br}{\operatorname{Br}}
\newcommand{\gal}{\operatorname{Gal}}
\newcommand{\tr}{\operatorname{Tr}}

\DeclareFontFamily{U}{mathx}{\hyphenchar\font45}
\DeclareFontShape{U}{mathx}{m}{n}{
      <5> <6> <7> <8> <9> <10>
      <10.95> <12> <14.4> <17.28> <20.74> <24.88>
      mathx10
      }{}
\DeclareSymbolFont{mathx}{U}{mathx}{m}{n}
\DeclareFontSubstitution{U}{mathx}{m}{n}
\DeclareMathAccent{\widecheck}{0}{mathx}{"71}
\DeclareMathAccent{\wideparen}{0}{mathx}{"75}

\renewcommand{\epsilon}{\varepsilon}

\newcommand{\cha}{\operatorname{char}}

\newcommand{\upic}{\underline{\operatorname{Pic}}}

\newcommand{\et}{\mathrm{\acute{e}t}}

\renewcommand{\tor}{\mathrm{tor}}

\newcommand{\ptor}{\Pi^{\rm tor}}
\newcommand{\pet}{\Pi^{\mathrm{\acute{e}t}}}
\newcommand{\unsim}{\mathord{\sim}}

\newcommand{\cact}{\operatorname{Cact}}
\newcommand{\vect}{\operatorname{Vect}}
\newcommand{\fet}{\operatorname{F\acute{e}t}}
\newcommand{\spf}{\operatorname{Spf}}

\begin{document}

\title[The section conjecture for the toric fundamental group over $p$-adic fields]{The section conjecture for the toric fundamental group over $p$-adic fields}

\begin{abstract}
	Loosely speaking, the toric fundamental group is the Tannaka dual of a category of vector bundles which become direct sums of line bundles on a finite étale cover. In characteristic $0$, it is an extension of the étale fundamental group scheme by a projective limit of tori.
	
	We prove the analogue of Grothendieck's section conjecture for the toric fundamental group over $p$-adic fields. As a consequence, we give a new interpretation of Selmer sections as toric Galois sections.
\end{abstract}

\author[Bresciani]{Giulio Bresciani}
\address{Department of Mathematics, Universit\`a di Pisa}
\email{giulio.bresciani1@unipi.it}
\thanks{The author was partially supported by the FIS grant ``RAPTbyFUN'' FIS-2024-04657, and by the PRIN project ``Derived and underived algebraic stacks and applications''.}


\maketitle
\section{Introduction}\label{sect:int}

Curves are smooth, projective and geometrically connected. A curve is hyperbolic if it has genus $\ge 2$.

\subsection{The section conjecture}

Given a geometrically connected scheme $X$ over a field $k$, there is a short exact sequence of étale fundamental groups
\[1 \to \pi_{1}(X_{k^{s}}) \to \pi_{1}(X) \to \gal(k^{s}/k) \to 1\]
where $k^{s}$ is a separable closure of $k$, and $\gal(k^{s}/k)=\pi_{1}(\spec k)$.

The space of Galois sections $\cS^{\et}_{X/k}$ is the set of sections of this exact sequence modulo the action by conjugation of $\pi_{1}(X_{k^{s}})$. The functoriality of $\pi_{1}$ induces a map
\[X(k) \to \cS^{\et}_{X/k},\]
sometimes called \emph{profinite Kummer map}. Galois sections in the image of this map are called \emph{geometric}.

In 1983, in a famous letter to Faltings \cite{grothendieck-faltings}, Grothendieck made the following conjecture, known as the \emph{section conjecture}.

\begin{conjecture}[Grothendieck]
	If $k$ is finitely generated over $\QQ$ and $X$ is a curve of genus $\ge 2$ over $k$, then $X(k)\to\cS^{\et}_{X/k}$ is bijective.
\end{conjecture}

Injectivity was already known to Grothendieck and follows from the Mordell--Weil theorem, surjectivity is still wide open.

\subsection{The toric section conjecture}

It is ``folklore knowledge'' that an étale fundamental group scheme (or Nori-étale fundamental group) can be defined as the Tannaka dual of the category of vector bundles trivialized by a finite étale cover (see \cite{nori-representations} for the analogous construction of Nori's fundamental group), and that the étale fundamental group scheme is essentially a rephrasing of Grothendieck's étale fundamental group in a different language. 

Let $X$ be a smooth, projective, geometrically connected variety over $k$, with an ample line bundle $H$. We define the toric fundamental group scheme as the Tannaka dual of the category of vector bundles $V$ with a splitting $V = \bigoplus_{q \in \QQ}V_{q}$ such that, for some connected finite étale cover $f: Y \to X$, the pullback $f^{*}V_{q}$ is a direct sum of line bundles of $H$-degree equal to $\deg(f) \cdot q$; the splitting is part of the datum, see \S\ref{sect:cactus}. An alternative construction can be given in greater generality, and it shows that the toric fundamental group does not depend on the choice of $H$, see \S\ref{sect:toricgerbe}.

Analogously to the étale case, there is a space $\cS^{\tor}_{X/k}$ of toric Galois sections and a natural map $X(k) \to \cS^{\tor}_{X/k}$, which we call the \emph{toric Kummer map}.

\begin{conjecture}\label{conj:torQ}
	Let $k$ be a finitely generated extension of $\QQ$. Let $X$ be a hyperbolic curve over $k$. The toric Kummer map
	\[X(k) \to \cS^{\tor}_{X/k}\]
	is bijective.
\end{conjecture}

We prove that, for curves and abelian varieties in characteristic $0$, $\cS^{\tor}_{X/k} \subseteq \cS^{\et}_{X/k}$, see Theorem~\ref{thm:toretinj}. As a consequence, Grothendieck's conjecture implies Conjecture~\ref{conj:torQ}.

\subsection{Main results}

We prove that the toric version of the conjecture holds over $p$-adic fields both for abelian varieties and hyperbolic curves.

\begin{theorem}\label{thm:torab}
	Let $A$ be an abelian variety $A$ over a field $k$ finite over $\QQ_{p}$. Then
	\[A(k) \to \cS^{\tor}_{A/k}\]
	is bijective.
\end{theorem}

\begin{theorem}\label{thm:torQp}
	Let $X$ be a smooth, projective hyperbolic curve over a field $k$ finite over $\QQ_{p}$. Then
	\[X(k) \to \cS^{\tor}_{X/k}\]
	is bijective.
\end{theorem}

The index of a curve is the greatest common divisor of the degrees of its closed points. Recall that J. Stix has proved that, if $\cS^{\et}_{X/k} \neq \emptyset$, then the index is a power of $p$ \cite{stix-period}.

We prove Theorem~\ref{thm:torQp} by reducing it to a statement about the original version of the conjecture, and then proving this statement. Here it is.

\begin{theorem}\label{thm:nbd-p}
	Let $X$ be a smooth, projective hyperbolic curve over a field $k$ finite over $\QQ_{p}$. If $s \in \cS^{\et}_{X/k}$ is a Galois section not associated with a $k$-point of $X$, there exists a finite étale cover $Y \to X$ with a lift $r \in \cS^{\et}_{Y/k}$ of $s$ such that $p$ divides the index of $Y$.
	
	If $X(k) = \emptyset$, we can choose $Y \to X$ as a geometrically Galois cover whose degree is a power of $p$.
\end{theorem}

More generally, in Theorem~\ref{thm:nbd-p} we can assume that $Y \to X$ is geometrically Galois of degree a power of $p$ if the pro-$p$ section induced by $s$ is not geometric, see \S\ref{sect:curves} for details. We recall that Y. Hoshi has proved the existence of non-geometric pro-$p$ sections \cite{hoshi-prop}.

It's not known whether Grothendieck's conjecture holds over $p$-adic fields. By a well-known result of A. Tamagawa \cite[Proposition 2.8]{tamagawa}, it is enough to study the case of curves without rational points. Theorem~\ref{thm:nbd-p} takes Tamagawa's result one step further, telling us that over $p$-adic fields we can actually focus on curves of index different from $1$.

While Grothendieck's conjecture is not known over $p$-adic fields, we mention that F. Pop and J. Stix proved a valuation-theoretic analogue \cite{pop-stix}.

Recall that a Galois section $s\in\cS^{\et}_{X/k}$ of a curve $X$ over a number field $k$ is Selmer if its base change to $k_{\nu}$ is geometric for every finite place $\nu$. Write $\operatorname{Sel}_{X/k}\subset\cS^{\et}_{X/k}$ for the subset of Selmer Galois sections. As a consequence of Theorem~\ref{thm:torQp}, we prove that toric and Selmer sections coincide.

\begin{theorem}\label{thm:selmer}
	Let $X$ be a smooth, projective curve of genus $\ge 2$ over a number field $k$. Then 
	\[\operatorname{Sel}_{X/k} = \cS_{X/k}^{\tor} \subseteq \cS^{\et}_{X/k}.\]
\end{theorem}

We remark that the toric fundamental group does \emph{not} behave well under base change of the base field, if the extension is not algebraic. Since $k_{\nu}/k$ is not algebraic, the proof of Theorem~\ref{thm:selmer} is less straightforward than it might seem at first glance.

We recall that A. Betts and J. Stix have recently proved that, if $k$ contains no CM subfields and $\nu$ is a finite place of $k$, the natural map $\operatorname{Sel}_{X/k}\to X(k_{\nu})$ has finite image \cite{betts-stix}.

\subsection{An unexpected byproduct}

While studying the toric section conjecture, we stumbled upon a purely geometric result which, up to our knowledge, was not previously known. We record it here, just to avoid it being buried in the main text.

\begin{theorem}\label{thm:etale-sep}
	Let $X$ be a smooth, projective curve of genus $\ge 2$ over an algebraically closed field of characteristic $0$. For every pair of distinct points $x_{1},x_{2} \in X$ and for every integer $n \ge 2$, there exists a finite étale cover $Y \to X$ with the following property.
	
	For every pair of points $y_{1}, y_{2} \in Y$ over $x_{1},x_{2}$ and every $m$ with $\gcd(n,m) = 1$, we have $m([y_{1}] - [y_{2}]) \neq 0$ in $\pic^{0}(Y)$.
\end{theorem}

Quite surprisingly, even if Theorem~\ref{thm:etale-sep} is basically a statement over $\CC$, our proof uses anabelian methods over $p$-adic fields. It is natural to ask whether Theorem~\ref{thm:etale-sep} is true in positive characteristic as well, and if one can prove more strongly that $[y_{1}]-[y_{2}]$ is not torsion (i.e. the case $n = 1$). 

\subsection{Acknowledgements} 
The original idea of this paper goes back to 2018: while reading Stix' paper \cite{stix-period}, I realized that some version of Theorem~\ref{thm:nbd-p} should have been true. During this time, I have talked about it to several mathematicians. I want to thank D. Abramovich, F. Binda, S. Brochard, L. Capuano, F. Gallinaro, Y. Hoshi, A. Javanpeykar, D. Lombardo, J. Stix, A. Tamagawa, A. Vezzani, A. Vistoli and L. Zhang for having the patience to hear my half-baked thoughts and for suggesting various approaches and references.  Special thanks to B. Poonen for pointing out to me the Tate--Voloch conjecture and its proof by T. Scanlon, which is an essential ingredient of this paper. I also want to mention that D. Lombardo and A. Tamagawa provided the proofs of a couple of propositions contained in previous versions of the paper; while these didn't make it to the final cut, they certainly pointed me in the right direction.

Finally, let me mention that a previous version of the paper contained a serious mathematical mistake. I want to give my most heartfelt thanks to Y. Hoshi and A. Tamagawa for finding it.

\subsection{AI disclosure}

The results of section~\ref{sect:tower} were obtained in collaboration with ChatGPT 5.6 (see Remark~\ref{rem:ai} for details). The first draft of section~\ref{sect:cactus} was written by ChatGPT under the guidance of the author, who informally described its contents, and then edited by the author. Sections~\ref{sect:int}-\ref{sect:nf} were written by the author. Finally, ChatGPT reviewed the paper and spotted several minor mistakes. The author takes full responsibility for the contents of the paper.

\subsection{Fundamental gerbes} N. Borne and A. Vistoli \cite{borne-vistoli-nori} \cite{borne-vistoli} developed an alternative approach to fundamental groups in algebraic geometry, i.e. \emph{fundamental gerbes}. The languages of fundamental groups and fundamental gerbes are equivalent, and translations can be made back and forth. Fundamental gerbes have the disadvantage of being more technical since they are stacks, but have the advantage of not requiring the choice of a base point.

Fundamental gerbes also have another advantage which is very important for us: Borne and Vistoli constructed a general framework which automatically produces a wealth of variants of the original étale fundamental gerbe, the toric fundamental gerbe being one of them. Conversely, toric fundamental groups are currently not present in the literature.

We use Borne and Vistoli's language in order to have a reference to rely on. See \S\ref{sect:intgerbes} for a short and informal introduction to the language of fundamental gerbes.

In \S\ref{sect:cactus}, we give a Tannakian interpretation of toric fundamental groups and gerbes. 

\subsection{On cohomology}

We frequently use cohomology of group schemes which are not of finite type. In order to do this, we work in a fixed Grothendieck universe and use the fpqc topology. If $G$ is a (possibly non-abelian) group scheme over $k$, $\H^{1}(k,G)$ is the set of fpqc torsors. If $G$ is abelian, $\H^{1}(k,G)$ has a natural structure of abelian group, and $\H^{2}(k,G)$ is the group of fpqc gerbes banded by $G$, see \S\ref{sect:abgerbes}

\subsection{Notation and conventions} Curves are smooth, projective and geometrically connected. A variety over a field $k$ is a separated, geometrically integral scheme of finite type over $k$. An algebraic group over $k$ is a group scheme of finite type over $k$.

For a proper, geometrically connected and geometrically reduced scheme $X$ over a field $k$, we write $\upic_{X}$ for the Picard scheme and $\pic(X)$ for the group of line bundles on $X$. There is a natural inclusion $\pic(X)\subseteq\upic_{X}(k)$, which is an equality if $X(k)\neq\emptyset$.

Given an abelian variety $A$ and a positive integer $n$, $A[\infty]$ is the subgroup of all torsion points, $A[n]$ is the kernel of multiplication by $n$, $A[n^{\infty}]$ are the torsion points whose order divides a power of $n$, and $A[n']$ are the torsion points whose order is prime with $n$.

Given a variety $X$ over $k$, the degree of a closed point $x \in X$ is the degree of the residue field $k(x)$ over $k$. If $k$ is a local field, the ramification of $x$ is the ramification index of $k(x)$ over $k$.

\section{Gerbes}\label{sect:gerbes}

Recall that, given a site $\cC$, a \emph{gerbe} $\cG$ on $\cC$ is a stack in groupoids such that
\begin{itemize}
	\item for every $S \in \cC$, any two objects of $\cG(S)$ are locally isomorphic, and
	\item for every $S \in \cC$, there exists a covering $\{S_{i} \to S\}$ such that $\cG(S_{i})$ is non-empty for every $i$.
\end{itemize}

An algebraic gerbe over a field $k$ is an algebraic stack over $k$ which is a gerbe for the fppf topology. We need to use gerbes which are projective limits of algebraic gerbes, and these are in general not gerbes for the fppf topology, but only for the fpqc one. Because of this, when we work with a gerbe over a field we always mean that it is a gerbe for the fpqc topology. Furthermore, we work with \emph{affine} gerbes over a field, i.e. fpqc gerbes with affine diagonal, possessing an affine chart.

We are going to use the concept of \emph{locally full} morphism of gerbes, and the \emph{canonical factorization} of a morphism of gerbes, which were introduced in \cite{borne-vistoli}. These are simply the gerbe-theoretic analogues of a surjective homomorphism of groups, and of the factorization of a homomorphism into a surjective homomorphism composed with an injective one.

\begin{definition}[{\cite[Definition 3.4]{borne-vistoli}}]
	Let $f: \Upsilon \to \Phi$ be a morphism of affine gerbes over a field $k$. We say that $f$ is \emph{locally full} if, for any extension $k'$ of $k$ and any object $\xi \in \Upsilon(k')$, the induced homomorphism of group schemes $\aut_{\Upsilon}(\xi) \to \aut_{\Phi}(f(\xi))$ is faithfully flat.
	
	Equivalently, $f$ is locally full if $\Upsilon$ is a relative gerbe over $\Phi$ \cite[Proposition 3.10]{borne-vistoli}.
\end{definition}

\begin{proposition}[{\cite[Proposition 3.9]{borne-vistoli}}]\label{prop:canonical}
	Let $f: \Upsilon \to \Phi$ be a morphism of affine gerbes over a field $k$. There exists a factorization $\Upsilon \to \Phi' \to \Phi$ with $\Upsilon \to \Phi'$ locally full and $\Phi' \to \Phi$ faithful.
	
	Furthermore, the factorization is canonical: if $\Upsilon \to \Phi'' \to \Phi$ is another such factorization, there exists an equivalence $\Phi' \to \Phi''$ and a $2$-commutative diagram
	\[\begin{tikzcd}[column sep=small, row sep=small]
														&	\Phi'\ar[dr]\ar[dd,dashed,crossing over]	&						\\
		\Upsilon\ar[ur]\ar[dr]							&								&	\Phi				\\
														&	\Phi''\ar[ur]				&	
	\end{tikzcd}\]\qed
\end{proposition}

\subsection{Abelian gerbes}\label{sect:abgerbes}

A gerbe $\cG$ on a site $\cC$ is \emph{abelian} if, for every $S \in \cC$ and every $s \in \cG(S)$, the automorphism group $\aut_{\cG}(s)$ is abelian.

If $\cG$ is an abelian gerbe over a small site, there exists a sheaf of groups $G$ called the \emph{band} of $\cG$ such that $\aut_{\cG}(s) \simeq G(S)$, plus some obvious compatibility properties \cite[\href{https://stacks.math.columbia.edu/tag/0CJY}{Tag 0CJY}]{stacks-project}.

A morphism $\cG \to \cH$ of abelian gerbes induces a homomorphism of bands $G \to H$.

The following is well-known, but the only reference we know works in the much more complex non-abelian case \cite{giraud}. We give an easy proof in the abelian case.

\begin{lemma}\label{lem:induced}
	Let $\cG$ be an abelian gerbe with band $G$ on a small site $\cC$, $H$ a sheaf of abelian groups and $f: G \to H$ a homomorphism of sheaves of groups.
	
	There exists a gerbe $\cG \times^{G}H$ banded by $H$, called the \emph{induced gerbe}, with a morphism $\cG \to \cG \times^{G}H$ whose homomorphism of bands is $f$. 
	
	Furthermore, if $\cG \to \cH$ is another such morphism, there exists an equivalence $\cG \times^{G}H \to \cH$ making the obvious diagram $2$-commutative and whose associated morphism of bands is the identity of $H$.
\end{lemma}

\begin{proof}
	Define a category fibered in groupoids $\cH_{0}$ whose objects coincide with the objects of $\cG$, and arrows defined as follows. If $s_{1} \in \cH_{0}(S_{1}) = \cG(S_{1})$, $s_{2}\in \cH_{0}(S_{2}) = \cG(S_{2})$ and $S_{1} \to S_{2}$ is a morphism in $\cC$, then
	\[\hom_{\cH_{0}}(s_{1}, s_{2}) \eqdef \hom_{\cG}(s_{1}, s_{2}) \times^{G}H = \left(\hom_{\cG}(s_{1}, s_{2}) \times H(S_{1})\right)/G(S_{1})\]
	where, for the fixed arrow $S_1\to S_2$, the contracted action is
\[
(u,h)\cdot g=(u\circ g,f(g)^{-1}h),
\qquad g\in G(S_1).
\]
Sending an arrow $u$ to the class of $(u,1)$ defines a morphism $\cG\to\cH_0$.
	
	Since $\cC$ is a small site, we may define $\cG \times^{G}H$ as the stackification of $\cH_{0}$ \cite[\href{https://stacks.math.columbia.edu/tag/02ZN}{Tag 02ZN}]{stacks-project}; we get a morphism $\cG \to \cH_{0} \to \cG \times^{G}H$ by composition.
	
	Given an object $S \in \cC$, any element of $\cG \times^{G}H(S)$ comes locally from $\cH_{0}$, and hence from $\cG$: since $\cG$ is a gerbe, this implies that any two elements of $\cG \times^{G}H(S)$ are locally isomorphic, proving that $\cG \times^{G}H$ is a gerbe. Let us show that it is abelian banded by $H$.
	
	For every element $s \in \cG \times^{G}H(S)$, there is an induced homomorphism of sheaves of groups $H \to \underaut_{\cG\times^G H}(s)$; we want to show that this is an isomorphism. Every element of $\underaut_{\cG\times^G H}(s)$ comes locally from an automorphism defined in $\cH_{0}$; this implies that $H \to \underaut_{\cG\times^G H}(s)$ is surjective. Let us prove injectivity.
	
	Up to passing to a cover, we may assume that $s$ lifts to $\cG(S)$, and up to replacing $\cC$ with the comma category $\cC/S$ we may also assume that $\cG$ is the classifying stack $\cB G$. Under this assumption, we have a natural morphism $\cH_{0} \to \cB H$ where $\cB H$ is the classifying stack of $H$, and hence an induced morphism $\cG \times^{G}H \to \cB H$. This gives a retraction $\underaut_{\cG\times^G H}(s) \to H$ of $H \to \underaut_{\cG\times^G H}(s)$, proving injectivity.
	
	If $\cG \to \cH$ is another morphism of gerbes whose morphism of bands is $f$, then there is a natural morphism $\cH_{0} \to \cH$ inducing a morphism $\cG \times^{G}H \to \cH$ whose morphism of bands is the identity of $H$. A morphism of abelian gerbes is an equivalence if and only if the corresponding morphism of bands is an isomorphism, hence we conclude.
\end{proof}

Given a sheaf of abelian groups $G$ on a small site $\cC$ and two gerbes $\cG$, $\cG'$ banded by $G$, we can define the sum $\cG \oplus \cG'$ as the induced gerbe $(\cG\times\cG')\times^{G \times G}G$ along the multiplication homomorphism $G\times G \to G$. Similarly, we can define an inverse gerbe $\cG^{-1}$ induced by the homomorphism $g \mapsto g^{-1}$.

This gives a group structure to the set of equivalence classes of gerbes banded by $G$; this group is naturally isomorphic to the sheaf cohomology group $\H^{2}(G)$ \cite[IV, \S 3.4, Theorem 3.4.2(i)]{giraud}.

We need to work with fpqc gerbes, but unfortunately the fpqc site is not small. Hence, we work in a fixed Grothendieck universe; this will allow us to apply Lemma~\ref{lem:induced} to fpqc gerbes.

\subsection{Projective limits of gerbes}

Consider a pro-algebraic fpqc gerbe $\Pi$ over $k$, i.e. a limit in the sense of \cite[\S 3]{borne-vistoli-nori} of affine, algebraic gerbes. By definition, $\Pi=\projlim \Gamma$ is the limit along a cofiltered poset $I$ of a pseudofunctor $\Gamma:I \to (\operatorname{AffGer}/k)$ cf. \cite[Remark 3.4, Definition 3.5]{borne-vistoli-nori}. This means that, for every $j \ge i \ge h$, there exists a morphism $\gamma_{i,j}:\Gamma_{j} \to \Gamma_{i}$ and a $2$ isomorphism $\lambda_{h,i,j}:\gamma_{h,i} \circ \gamma_{i,j} \Rightarrow \gamma_{h,j}$. Furthermore, these satisfy a cocycle condition.

\begin{definition}
	Given a pro-algebraic gerbe $\Pi=\projlim \Gamma$, we define the set $\Pi(k)_{\rmq}$ of \emph{quasi-sections} as the projective limit of the sets $\Gamma_{i}(k)/\unsim$ of isomorphism classes of sections; there is an obvious map $\Pi(k)/\unsim \to \Pi(k)_{\rmq}$.

It is straightforward to check that this definition does not depend on the choice of the presentation.
\end{definition}

The following is the stacky analogue of \cite[Lemma 44]{stix}.

\begin{lemma}\label{lem:profcont}
	Assume that $\Pi=\projlim \Gamma$ is profinite, i.e. $\Gamma_{i}$ is a finite affine gerbe in the sense of \cite[\S 4]{borne-vistoli-nori} for every $i$. Then $\Pi(k)/\unsim\to\Pi(k)_{\rmq}$ is bijective.
\end{lemma}

\begin{proof}
	If $s\in\Pi(k)$, then $\Pi(k)/\unsim\simeq\H^{1}(k,\aut_{\Pi}(s))$ and the statement is proved in \cite[Proposition 2.3]{bresciani-proed}, hence injectivity is automatic. We have to prove surjectivity.

	Let $(s_{i})_{i}\in \Pi(k)_{\rmq}$ be a quasi-section; by definition, for every $j \ge i$ there exists an isomorphism $\sigma_{i,j}:\gamma_{i,j}(s_{j})\to s_{i}$. If these isomorphisms satisfy the cocycle condition
	\[\sigma_{h,j} = \sigma_{h,i} \circ \gamma_{h,i}(\sigma_{i,j}) \circ \lambda_{h,i,j}(s_{j})^{-1} : \]
	\[ : \gamma_{h,j}(s_{j}) \to \gamma_{h,i} \circ \gamma_{i,j}(s_{j}) \to \gamma_{h,i}(s_{i}) \to s_{h}\]
	then we get the desired section $s\in\Pi(k)$ by definition of projective limit; the cocycle condition is in general not satisfied, so we have to modify the isomorphisms $\sigma_{i,j}$.
	
	Denote by $\rF\rP(I)$ the poset of finite sub-posets of $I$ ordered by inclusion, it is cofiltered as well. For every $S\in\rF\rP(I)$, denote by $\Lambda(S)$ the set of coherent choices of isomorphisms $\sigma_{i,j}':\gamma_{i,j}(s_{j})\to s_{i}$ for $j \ge i$ in $S$, i.e. an element of $\Lambda(S)$ is a choice of $\sigma_{i,j}'$ for every $j \ge i\in S$ such that the $\sigma_{i,j}'$s satisfy the cocycle condition written above restricted to the elements of $S$.
	
	If $S\subset S'$, we have an obvious forgetful map $\Lambda(S') \to \Lambda (S)$. Since $S$ is finite and, by hypothesis, there is a finite number of isomorphisms $\gamma_{i,j}(s_{j})\to s_{i}$, we get that $\Lambda(S)$ is finite. I claim that $\Lambda(S)$ is non-empty for every $S$.
	
	In fact, since $I$ is cofiltered up to enlarging $S$ we can assume that it has a maximum $m$. For every $j \ge i$ in $S$, define
	\[\sigma'_{i,j} \eqdef \sigma_{i,m} \circ \lambda_{i,j,m}(s_{m}) \circ \gamma_{i,j}(\sigma_{j,m}^{-1}):\]
	\[: \gamma_{i,j}(s_{j}) \to \gamma_{i,j}\circ \gamma_{j,m}(s_{m}) \to \gamma_{i,m}(s_{m}) \to s_{i}.\]
	It is straightforward to check that these satisfy the cocycle condition in $S$ and thus define an element of $\Lambda(S)$.
	
	The projective limit of $\Lambda$ along $\rF\rP(I)$ is the set of coherent choices of isomorphisms $\gamma_{i,j}(s_{j})\to s_{i}$ for every $j \ge i$ in $I$. A projective limit of finite, non-empty sets is non-empty, hence such a coherent choice always exists.
\end{proof}

\section{Fundamental gerbes}

In this section, we review and integrate some material from \cite{borne-vistoli-nori} \cite{borne-vistoli}.

\subsection{An informal introduction to fundamental gerbes}\label{sect:intgerbes}

We start by giving a short and informal introduction to fundamental gerbes. We refer to \cite{borne-vistoli-nori, borne-vistoli} for precise statements and proofs.

Let $\operatorname{HoTop}$ be the homotopy category of connected CW complexes. Given $X\in\operatorname{HoTop}$ and a base point $x$, its fundamental group $\pi_1(X,x)$ can be defined using based loops up to homotopy, or as the group of deck transformations of the universal cover.

One can also characterize the classifying space $\cB\pi_1(X,x)$ directly, without first choosing a base point. Let $\operatorname{EM}_1\subseteq\operatorname{HoTop}$ be the full subcategory of Eilenberg--MacLane spaces in degree $1$, namely connected CW complexes whose higher homotopy groups vanish. There is a left adjoint $\Pi:\operatorname{HoTop}\to\operatorname{EM}_1$, characterized by the natural bijections
\[
\hom(X,Y)=\hom(\Pi_X,Y),\qquad Y\in\operatorname{EM}_1,
\]
in the homotopy category.

The space $\Pi_X$ can be constructed as $\cB\pi_1(X,x)$, or by attaching cells to $X$ to kill its higher homotopy groups \cite[Example 4.17]{hatcher}. Its universal property determines it up to unique isomorphism in the homotopy category, independently of the construction and of the base point.

The idea of fundamental gerbes is to repeat this construction in the algebro-geometric setting. There is a problem, though: what is the right analogue of Eilenberg--MacLane spaces in degree $1$? We cannot expect them to be schemes, even if we work over $\CC$: two elliptic curves are homotopically equivalent, but there is no simple way to see this in the category of algebraic varieties. Algebraic varieties are just too ``rigid'' to do homotopy.

A better choice is to use stacks, and in particular (pro-)algebraic \emph{gerbes}. For instance, over $\CC$ every algebraic gerbe is a quotient stack
\[[\spec \CC / G]\]
for some algebraic group $G$ acting trivially on $\spec \CC$, and this quotient is by definition the classifying stack $\cB G$.

Over an arbitrary base field $k$, an algebraic gerbe $\Phi$ is a classifying stack if and only if $\Phi(k)\neq\emptyset$. If $p\in\Phi(k)$, then $\aut_{\Phi}(p)$ is an algebraic group and $\Phi\simeq \cB\aut_{\Phi}(p)=[\spec k/\aut_{\Phi}(p)]$. Every algebraic gerbe over $k$ has a geometric point defined over the algebraic closure $\bar{k}$, hence we might think of gerbes as twisted forms of classifying stacks.

Notice that $\spec \CC$, seen as a topological space, is \emph{contractible}: in topology, the standard construction of the classifying space $\cB G$ is to take the quotient of a contractible topological space by a free action of $G$. With topological spaces, we need the action to be free, otherwise the quotient just doesn't have the desired properties.

With stacks, we do not need the action to be free: quotient stacks are defined in such a way that we can \emph{pretend} that the action is always free. In fact, one could argue that one of the main reasons for developing the theory of stacks is exactly to pretend that some non-free actions are free. Hence, the most natural choice to define classifying stacks is to take the simplest action on the simplest variety, i.e. the trivial action on one point.

So we replace Eilenberg--MacLane spaces with gerbes. But \emph{which} gerbes? This depends on the type of fundamental group we want to construct. The natural thing to do, rather than choosing the gerbes themselves, is to choose which automorphism groups of points of gerbes we allow. This basically amounts to choosing for which groups $G$ we want to consider the classifying stack $\cB G$.

The first and most natural choice is to consider finite étale group schemes: this will lead us to defining the étale fundamental gerbe $X \to \Pi^{\et}_{X/k}$, which is a close relative of the étale fundamental group. The other class of groups we are interested in is virtually toric groups, i.e. algebraic groups whose connected component is of multiplicative type, this will lead us to define the toric fundamental gerbe $X \to \Pi_{X/k}^{\tor}$.

There are many other choices, though. If we consider finite group schemes, we get Nori's fundamental gerbe, the analogue of Nori's fundamental group. Other choices are abelian algebraic groups, or unipotent algebraic groups, see \cite{borne-vistoli}.

Let $\cT$ be such a class of group schemes, and let $\Pi^{\cT}_{X/k}$ be the associated fundamental gerbe (there are some assumptions to make on $\cT$ in order to have a functioning machinery, see \cite{borne-vistoli}). 

If $X$ has a \emph{rational} base point $x\in X(k)$, then we can repeat a construction of Nori \cite[Chapter II]{nori} and define a $\cT$-group scheme $\pi_{1}^{\cT}(X,x)$: the fundamental gerbe $\Pi_{X/k}^{\cT}$ coincides with the classifying stack $\cB\pi_{1}^{\cT}(X,x)$. There is a universal $\pi_{1}^{\cT}(X,x)$-torsor $U\to X$, the set of isomorphism classes of $\Pi^{\cT}_{X/k}(k)$ coincides with $\H^{1}(k,\pi_{1}^{\cT}(X,x))$ and the morphism $X \to \Pi^{\cT}_{X/k}$ maps a point $p\in X(k)$ to the fiber $U_{p}$, which is a torsor over $\spec k$ and corresponds to a class in $\H^{1}(k,\pi_{1}^{\cT}(X,x))=\Pi^{\cT}_{X/k}(k)/\unsim$.

The point of using gerbes is exactly that we do not need any base point at all, not even a geometric one. 

Moreover, the section conjecture has a particularly simple formulation in the language of fundamental gerbes. If $\Pi^{\et}_{X/k}$ is the étale fundamental gerbe, the space of Galois sections is simply the set of isomorphism classes of $\Pi^{\et}_{X/k}(k)$, and the profinite Kummer map is just the evaluation of $X \to \Pi^{\et}_{X/k}$ on $k$-points.

Finally, let us mention that the fundamental gerbe will typically not exist in the class of \emph{algebraic} gerbes, for the same reason why we cannot expect the étale fundamental group to exist as a finite group. There is a limit argument to make: just as the étale fundamental group is profinite, fundamental gerbes are pro-algebraic.

\subsection{The étale fundamental gerbe}

Given a geometrically connected fibered category $X$ over $k$ cf. \cite[Appendix]{bresciani-anab}, the étale fundamental gerbe is a profinite étale gerbe $\pet_{X/k}$ with a morphism $X\to\pet_{X/k}$ which is universal among morphisms $X \to \Phi$ where $\Phi$ is a finite étale gerbe \cite[\S 8]{borne-vistoli-nori} \cite[Appendix]{bresciani-anab}.

When $X$ is a quasi-compact, quasi-separated scheme, the space $\cS^{\et}_{X/k}$ of Galois sections of $X$ is in natural bijection with $\pet_{X/k}(k)/\unsim$ \cite[Proposition 9.3]{borne-vistoli-nori} \cite[Proposition A.19]{bresciani-anab}. By Lemma~\ref{lem:profcont}, it is in natural bijection with $\pet_{X/k}(k)_{\rmq}$ as well.

Base change results for the étale fundamental gerbe are scattered through the literature. Let us collect them all in one place.

Recall that a fibered category is \emph{concentrated} if there exists an affine scheme $U$ with a representable, faithfully flat, quasi-compact, quasi-separated morphism $U \to X$.

\begin{theorem}\label{thm:et-bc}
	Let $X$ be a concentrated, geometrically connected fibered category over $k$, $k'/k$ a field extension, $X' = X_{k'}$. The base change morphism 
	\[\tau:\pet_{X'/k'} \to \pet_{X/k}\times_{k}k'\]
	is an isomorphism under any of the following assumptions.
	\begin{enumerate}
		\item $k'/k$ is algebraic.
		\item $\cha k = 0$.
		\item $X$ is a proper, geometrically connected scheme.
	\end{enumerate}
\end{theorem}

\begin{proof}
	Assume first that $k'/k$ is algebraic. If $k'/k$ is algebraic and separable, $\tau$ is an isomorphism by \cite[Proposition A.18]{bresciani-anab} \cite[Theorem 9.3]{borne-vistoli}. We may then assume that $k$ is separably closed, and that $k'/k$ is purely inseparable. 
	
	If $\Phi$ is a finite, étale gerbe and $k$ is separably closed, there is a unique $\bar{k}$-point of $\Phi$ up to isomorphism which descends to a unique $k$-point of $\Phi$ (this follows easily by descent theory using the fact that Isom-schemes for $\Phi$ are finite étale). Since $|\Phi(k)/\sim| = 1$ if $\Phi$ is finite étale, the same is true in the case $\Phi$ is pro-finite étale thanks to the fact that $\Phi(k)/\sim = \Phi(k)_{\rm q}$ by Lemma~\ref{lem:profcont}. 
	
	This implies that $\pet_{X/k} = B_{k}G$, $\pet_{X'/k'} = B_{k'}G'$ for some pro-finite étale group schemes $G$, $G'$ over $k$, $k'$. Since $k$, $k'$ are separably closed, the group schemes $G$, $G'$ are pro-finite constant.
		
	Since $k'/k$ is purely inseparable, pullback defines an equivalence $\fet_{X} \to \fet_{X'}$ of categories of finite étale covers of $X$ and $X'$. This is proved in \cite[\href{https://stacks.math.columbia.edu/tag/04DZ}{Tag 04DZ}]{stacks-project} for schemes; the statement for fibered categories follows directly. Pullback also defines equivalences $\fet_{B_{k}G} \to \fet_{X}$, $\fet_{B_{k'}G'} \to \fet_{X'}$. This implies that $\fet_{B_{k}G} \to \fet_{B_{k'}G'}$ is an equivalence as well. Since $\fet_{B_{k}G}$ corresponds to the category of finite sets with a continuous $G$-action, we get that $G' \simeq G_{k'}$ by Grothendieck's classical theory of Galois categories. It follows that $\tau$ is an isomorphism.
	
	If $\cha k = 0$, this is \cite[Proposition A.23]{bresciani-anab}.
	
	Assume now that $X$ is a proper, geometrically connected scheme. By the first case, we may reduce to the case in which $k'$, $k$ are both algebraically closed. Then $\pet_{X/k} = B_{k}\pi_{1}(X)$, $\pet_{X'/k'} = B_{k'}\pi_{1}(X')$. Since $\pi_{1}(X') \simeq \pi_{1}(X)$ by \cite[Exposé X, Corollaire 1.8]{sga1}, we conclude that $\tau$ is an isomorphism.
\end{proof}

\subsection{Firm fibered categories}

While the étale fundamental gerbe exists for geometrically connected fibered categories, other types of fundamental gerbes need stronger assumptions cf. \cite[Theorem 7.1]{borne-vistoli}. Let us give these assumptions a name.

\begin{definition}
	A fibered category $X$ over a field $k$ is \emph{firm} if it is concentrated, geometrically reduced cf. \cite[Definition 4.3]{borne-vistoli}, and $\H^{0}(X, \cO_{X}) = k$.
\end{definition}

For instance a proper, geometrically connected, geometrically reduced scheme is firm.

\begin{lemma}\label{lem:firmet}
	Let $X$ be a firm fibered category over a field $k$. A profinite étale cover $Y \to X$ is firm if and only if it is geometrically connected. In particular, $X$ is geometrically connected.
\end{lemma}

\begin{proof}
	If $Y$ is firm then $\H^{0}(Y, \cO_{Y}) = k$; in particular, $Y$ is connected. Moreover, $Y$ is quasi-compact and quasi-separated, hence global sections commute with base change and $Y$ is geometrically connected. 
	
	Assume that $Y$ is geometrically connected.	Since $Y$ is profinite étale over $X$, then $Y$ is quasi-compact, quasi-separated and geometrically reduced. We have to prove that $\H^{0}(Y, \cO_{Y}) = k$.
	
	Since $Y$ is profinite over $X$, then $\H^{0}(Y,\cO_{Y})$ is algebraic over $\H^{0}(X, \cO_{X}) = k$. Since $X$ is geometrically reduced and $f$ is profinite étale, then $Y$ is geometrically reduced, hence $\H^{0}(Y, \cO_{Y})$ is an algebraic, separable extension of $k$. Since $Y$ is geometrically connected, we conclude that $\H^{0}(Y,\cO_{Y}) = k$.
\end{proof}

\subsection{The toric fundamental gerbe}\label{sect:toricgerbe}

We define the toric fundamental gerbe over fields of arbitrary characteristic, but study its properties in characteristic $0$.

An affine group scheme $G$ of finite type over a field $k$ is virtually toric, vt for short, if it is an extension of a finite étale group scheme by a connected group scheme of multiplicative type. Equivalently, $G^0$ is of multiplicative type; in positive characteristic it need not be reduced.

An affine gerbe cf. \cite[\S 3]{borne-vistoli} $\Phi$ is vt if the automorphism groups of its geometric points are vt. Pro-vt groups (resp. gerbes) are projective limits of vt groups (resp. gerbes).

\begin{definition}
	The \emph{toric fundamental gerbe} of a firm fibered category $X$ over a field $k$ is a pro-vt gerbe $\ptor_{X/k}$ with a morphism $X \to \ptor_{X/k}$ which is universal among morphisms $X \to \Phi$ where $\Phi$ is a vt gerbe \cite[Theorem 7.1]{borne-vistoli}, i.e. it is the fundamental gerbe for the class of vt groups in the terminology of Borne and Vistoli.
\end{definition}

The toric fundamental gerbe behaves well with respect to algebraic extensions of the base field.

\begin{theorem}[{\cite[Theorem 9.3]{borne-vistoli}}]\label{thm:basechange}
	Let $X$ be a firm fibered category over a field $k$, and $k'/k$ an algebraic, separable extension. The natural morphism
	\[\ptor_{X_{k'}/k'} \to \ptor_{X/k} \times_{k} k'\]
	is an isomorphism.
\end{theorem}

The étale fundamental gerbe behaves well for arbitrary extensions of the base field in characteristic $0$ \cite[Proposition A.23]{bresciani-anab}, but this fails for the toric fundamental gerbe. 
 
\begin{example}\label{ex:badchange}
	Let $\tilde{E} \to E$ be the universal cover of an elliptic curve defined over $\bar{\QQ}$. Consider the projective system $(E)_{n \in \NN}$ of $\NN$ copies of $E$ where the $nm$-th copy of $E$ maps to the $n$-th one by multiplication by $m$, so that $\tilde{E} = \projlim_{n}E$. By \cite[Théorème 8.5.2]{ega-iv.3}, the Picard group $\pic(\tilde{E})$ is the inductive limit $\indlim_{n}\pic(E) = \QQ \oplus E(\bar{\QQ})/E[\infty]$, where $E[\infty] \subset E(\bar{\QQ})$ is the subgroup of torsion points. For the same reason, $\pic(\tilde{E}_{\CC}) = \QQ \oplus E(\CC)/E[\infty]$.
	
	As we will prove in Lemma~\ref{lem:pichom}, $\ptor_{\tilde{E}/\bar{\QQ}}$ and $\ptor_{\tilde{E}_{\CC}/\CC}$ are abelian banded by pro-algebraic groups $T_{\tilde{E}}$ over $\bar{\QQ}$ and $T_{\tilde{E}_{\CC}}$ over $\CC$ satisfying
	\[\hom(T_{\tilde{E}},\GG_{m,\bar{\QQ}}) \simeq \pic(\tilde{E}) = \QQ \oplus E(\bar{\QQ})/E[\infty],\]
	\[\hom(T_{\tilde{E}_{\CC}},\GG_{m,\CC}) \simeq \pic(\tilde{E}_{\CC}) = \QQ \oplus E(\CC)/E[\infty].\]
	The former is countable while the latter is not, hence $T_{\tilde{E}_{\CC}}$ is not isomorphic to $T_{\tilde{E}} \times_{\bar{\QQ}} \CC$ and $\ptor_{\tilde{E}_{\CC}/\bar{\CC}}$ is not isomorphic to $\ptor_{\tilde{E}/\bar{\QQ}} \times_{\bar{\QQ}} \CC$.
\end{example}

\begin{definition}\label{def:torsect}
	The space $\cS_{X/k}^{\tor}$ of toric Galois sections of a firm fibered category $X$ over a field $k$ is the set $\ptor_{X/k}(k)_{\rmq}$ of quasi-sections of the toric fundamental gerbe.
\end{definition}

For a firm fibered category $X$ over a field $k$, define its relative Brauer group by
\[
\br(X/k):=\ker\bigl(\br(k)\to\H^2_{\mathrm{fppf}}(X,\GG_m)\bigr).
\]

\begin{lemma}\label{lem:brauerobstruction}
	Let $X$ be a firm fibered category over a field $k$. If $\br(X/k)\neq 0$, then $\cS^{\tor}_{X/k}=\emptyset$.
\end{lemma}

\begin{proof}
	A nonzero class in $\br(X/k)$ defines a gerbe $\cG$ banded by $\GG_m$ with $\cG(k)=\emptyset$ and a morphism $X\to\cG$. By the universal property, this morphism factors through $\ptor_{X/k}$. Since $\cG$ is of finite type, the morphism $\ptor_{X/k}\to\cG$ factors through a quotient of finite type, so a toric quasi-section would give a $k$-point of $\cG$, a contradiction.
\end{proof}

\begin{remark}
	One might wonder why define $\cS_{X/k}^{\tor}$ as $\ptor_{X/k}(k)_{\rmq}$, and not as $\ptor_{X/k}(k)/\unsim$. In the étale case, by Lemma~\ref{lem:profcont} the two sets always coincide, hence the choice is irrelevant. However, in the pro-algebraic case the choice is relevant, because equality does not hold even in fairly simple cases like the following Example~\ref{ex:nolim}.
	
	Example~\ref{ex:nolim} also shows that the set of isomorphism classes of sections of a pro-algebraic gerbe can be quite pathological, whereas quasi-sections have better behaviour. Furthermore, under an hypothesis satisfied by curves of genus $\ge 1$, the map $\ptor_{X/k}(k)_{\rmq}\to \pet_{X/k}(k)_{\rmq}$ is injective, see Theorem~\ref{thm:toretinj}, whereas the fibers of $\ptor_{X/k}(k)/\unsim \to \pet_{X/k}(k)/\unsim$ can be huge.
\end{remark}

\begin{example}\label{ex:nolim}
	Fix a prime $\ell$ and a field $k$, write $\SS_{\ell}$ for the $\ell$-adic solenoid, i.e. the projective limit of
	\[\dots \xrightarrow{\bullet^{\ell}} \GG_{m} \xrightarrow{\bullet^{\ell}} \GG_{m} \xrightarrow{\bullet^{\ell}} \GG_{m}.\]
	Clearly, $\projlim\H^{1}(k,\GG_{m}) = \cB \SS_{\ell}(k)_{\rmq}$ is trivial; we want to show that $\H^{1}(k,\SS_{\ell}) = \cB \SS_{\ell}(k)/\unsim$ is in general non-trivial.
	
	Consider the projection $\SS_{\ell}\to\GG_{m}$ to the first group in the tower, the kernel is $\ZZ_{\ell}(1)$. We have a long exact sequence in cohomology
	\[\GG_{m}(k)=k^{*}\to \H^{1}(k,\ZZ_{\ell}(1)) \to \H^{1}(k,\SS_{\ell}) \to \H^{1}(k,\GG_{m})=0.\]
	By \cite[Proposition 2.3]{bresciani-proed}, we have
	\[\H^{1}(k,\ZZ_{\ell}(1))=\projlim_{i}\H^{1}(k,\mu_{\ell^{i}})=\projlim_{i}k^{*}/k^{*\ell^{i}}.\]
	It follows that $\H^{1}(k,\SS_{\ell})$ is trivial if and only if the natural map $k^{*} \to \projlim_{i}k^{*}/k^{*\ell^{i}}$ is surjective. This is often false.
	
	For instance, the map is not surjective if $k=\QQ_{p}$, or more generally when $k$ has a surjective valuation $v:k^{*} \to \ZZ$. If $\alpha\in\ZZ_{\ell}\setminus\ZZ$ is an $\ell$-adic number not in $\ZZ$, and $x \in k^{*}$ has valuation $1$, then $x^{\alpha}$ is a well defined element of $\projlim_{i} k^{*}/k^{*\ell^{i}}$ not coming from $k^{*}$. In fact, the valuation extends to a map $v:\projlim_{i} k^{*}/k^{*\ell^{i}} \to \ZZ_{\ell}$ and $v(x^{\alpha})=\alpha$, whereas $v(y) \in \ZZ$ for every $y \in k^{*}$.
\end{example}

By applying Lemma~\ref{lem:profcont}, we immediately get the following.

\begin{lemma}
	Let $X$ be a firm fibered category over a field $k$. The natural map $\ptor_{X/k}\to\pet_{X/k}$ induces a map $\cS^{\tor}_{X/k}\to\cS^{\et}_{X/k}$. \qed
\end{lemma}

As for classical Galois sections, it is important to understand the behaviour of the space of toric Galois sections with respect to finite étale covers.

\begin{proposition}\label{prop:cartesian}
	Let $X$ be a firm fibered category over a field $k$, and $Y \to X$ a geometrically connected, profinite étale cover. Both squares of the natural $2$-commutative diagram
	\[\begin{tikzcd}
		Y\rar\dar	&	\ptor_{Y/k}\rar\dar		&	\pet_{Y/k}\dar		\\
		X\rar		&	\ptor_{X/k}\rar			&	\pet_{X/k}
	\end{tikzcd}\]
	are $2$-cartesian.
\end{proposition}

Before giving the actual proof, let us sketch the basic idea. For simplicity, assume $k = \bar{k}$ and that $Y \to X$ is a Galois cover with Galois group $G$. The proposition reduces to proving that, for every torsor $T \to Y$ for a virtually toric group $H$, there exists a ``Galois closure'', i.e. a torsor over $X$ for a virtually toric group with an equivariant map to $T$. For every $g \in G$, the pullback $g^{*}T$ is an $H$-torsor, and the direct product $H^{G}$ of $|G|$ copies of $H$ has a natural action of $G$ which permutes the components. Then $\prod_{g \in G}g^{*} T$ has a natural action of $H^{G} \rtimes G$ turning it into a $H^{G} \rtimes G$-torsor over $X$, and a $H^{G}$-torsor over $Y$. The natural projection $\prod_{g \in G}g^{*} T \to T$ is equivariant with respect to the projection $H^{G} \to H$ to the first component. In the proof, we basically repeat this argument at the level of gerbes.

\begin{proof}
	Assume first that $Y \to X$ is finite Galois with Galois group $G$, we are going to prove that $\ptor_{X/k} = [\ptor_{Y/k} / G]$ (see \cite{romagny} for the theory of group actions on stacks). The action on $Y$ induces an action on $\ptor_{Y/k}$, let us first show that the quotient stack $[\ptor_{Y/k} / G]$ is a pro-vt gerbe.
	
	Let $\ptor_{Y/k} \to \Upsilon$ be a morphism, where $\Upsilon$ is a vt gerbe of finite type over $k$; we have that $\Upsilon^{G}=\Upsilon \times \dots \times \Upsilon$ is a vt gerbe with a natural action of $G$, and there is a natural $G$-equivariant morphism $\ptor_{Y/k} \to \Upsilon^{G}$. Let $\ptor_{Y/k} \to \Upsilon' \to \Upsilon^{G}$ be the canonical factorization as in Proposition~\ref{prop:canonical}. The unicity of this factorization induces an action of $G$ on $\Upsilon'$ such that $\ptor_{Y/k} \to \Upsilon'$ is $G$-equivariant.
	
	Hence, up to replacing $\Upsilon$ with $\Upsilon'$, we may assume that we have a $G$-action on $\Upsilon$ such that $\ptor_{Y/k} \to \Upsilon$ is $G$-equivariant, and that $\ptor_{Y/k} \to \Upsilon$ is locally full. This implies that $\ptor_{Y/k}$ can be written as a projective limit of such gerbes $\Upsilon$, and $[\ptor_{Y/k} / G]$ is the limit of the gerbes $[\Upsilon /G]$; this shows that $[\ptor_{Y/k} / G]$ is pro-vt.
	
	Since $Y \to \ptor_{Y/k}$ is $G$-equivariant, we have an induced morphism $X \to [\ptor_{Y/k} / G]$, and hence a morphism $\ptor_{X/k} \to [\ptor_{Y/k} / G]$ by the universal property of the toric fundamental gerbe. On the other hand, $\ptor_{Y/k} \to \ptor_{X/k}$ is invariant with respect to the action of $G$, hence we get a map $[\ptor_{Y/k} / G] \to \ptor_{X/k}$. It is straightforward to check that these maps are inverses, hence we may identify $\ptor_{X/k} = [\ptor_{Y/k} / G]$. With the same argument, $\pet_{X/k} = [\pet_{Y/k} / G]$. 
	
	We now drop the assumption that $Y \to X$ is Galois, and only assume that $Y \to X$ is finite étale. Write $Y' = X \times_{\ptor_{X/k}} \ptor_{Y/k}$ and $\Pi = \ptor_{X/k} \times_{\pet_{X/k}} \pet_{Y/k}$. We want to show that the induced maps $Y \to Y'$, $\ptor_{Y/k} \to \Pi$ are isomorphisms, and we may do so after base changing to a separable closure of $k$ thanks to Theorem~\ref{thm:basechange}.
	
	We may then assume that $k$ is separably closed. There exists a fibered category $Z$ with a finite étale morphism $Z \to Y$ such that $Z \to Y \to X$ is Galois: in fact, the degree of $Y \to X$ is well-defined \cite[Lemma A.24]{bresciani-anab}, and then we may choose $Z$ as a suitable connected component of $Y \times_{X} \dots \times_{X} Y$ (connected components of fibered categories are defined via idempotents of $\H^{0}(X, \cO_{X})$ \cite[\S A.1]{bresciani-anab}). Denote by $G$ the Galois group of $Z \to X$ and $H \subset G$ the Galois group of $Z \to Y$. Since $k$ is separably closed, $Z$ is geometrically connected, and hence firm.
	
	By the above, $\ptor_{Y/k} = [\ptor_{Z/k} / H]$ and $\ptor_{X/k} = [\ptor_{Z/k} / G]$, and similarly for the étale fundamental gerbe. Every square in the diagram
	\[\begin{tikzcd}[column sep=small, row sep=small]
		Z\rar\dar							&	\ptor_{Z/k}\rar\dar										&	\pet_{Z/k}\dar\rar						&	\spec k\dar	\\
		Y \rar\dar	&	\left[ \ptor_{Z/k} / H \right]\rar\dar		&	\left[ \pet_{Z/k} / H \right]\dar\rar	&	\cB H\dar	\\
		X \rar		&	\left[ \ptor_{Z/k} / G \right]\rar			&	\left[ \pet_{Z/k} / G \right]\rar		&	\cB G		\\
	\end{tikzcd}\]
	is easily checked to be $2$-cartesian, hence we conclude.
	
	Finally, assume that $Y = \projlim_{i} Y_{i}$ is profinite étale over $X$, with $Y_{i} \to X$ finite étale. It is straightforward to prove that $\ptor_{Y/k} = \projlim_{i}\ptor_{Y_{i}/k}$ and $\pet_{Y/k} = \projlim_{i}\pet_{Y_{i}/k}$, hence the profinite case follows from the finite one.
\end{proof}

In order to define a \emph{toric fundamental group}, one has to fix a \emph{rational} base point of $X$, not a geometric one. This amounts to the fact that we want to obtain a group scheme over $k$, rather than an extension of $\gal(k^{s}/k)$.

More generally, we can fix a rational point of the fundamental gerbe.

\begin{definition}
	Let $\cT$ be a type of fundamental gerbe (e.g. étale or toric), and $s: \spec F \to \Pi^{\cT}_{X/k}$ a section for some field extension $F/k$. The $\cT$-fundamental group $\pi_{1}^{\cT}(X,s)$ is the automorphism group scheme of $s$. 
\end{definition}

If $F$ is algebraically closed, the $\cT$-fundamental group does not depend on the choice of $s$ since $\Pi^{\cT}_{X/k}$ is a gerbe and any two geometric sections are isomorphic, hence we might simply write $\pi_{1}^{\cT}(X)$ or even $\pi_{1}^{\cT}$. Over arbitrary fields, $\cT$-fundamental groups based at different sections are twisted forms of each other.

\begin{remark}
	If one wants to obtain an object which is more directly comparable with the étale fundamental group, one can choose a $\bar{k}$-base point $x:\spec \bar{k} \to \ptor_{X/k}$, and consider the group of $2$-commutative diagrams
	\[\begin{tikzcd}[column sep=tiny]
		\spec \bar{k} \ar[rr]\ar[dr,swap,"x"]	&	&	\spec \bar{k}\ar[dl,"x"]	\\
							&	\ptor_{X/k}	&
	\end{tikzcd}\]
	where the group operation is given by juxtaposition of diagrams. The group of diagrams is an extension of the Galois group $\gal(\bar{k}/k)$ by a pro-vt group over $\bar{k}$. If we consider the group of diagrams of the étale fundamental gerbe, and $x$ is the section corresponding to a point $x \in X(\bar{k})$, then the group of diagrams coincides with the classical étale fundamental group with base point $x$.
\end{remark}

\section{Torsion free abelian group schemes}

The connected component of the toric fundamental group is often torsion free, and simply connected: even though a non-trivial torus always has torsion, in the limit process defining the toric fundamental group it can happen that we dominate each finite étale cover, thus obtaining something which is torsion free and simply connected.

Because of this, we first study affine, torsion free group schemes beyond the usual finite type assumption.
	
A relative abelian group scheme $G \to S$ is torsion free if the kernel of the multiplication by $n$ map is a trivial relative group scheme for every $n$.

Recall that a homomorphism of affine group schemes $G \to Q$ over $k$ is a quotient if the corresponding map of Hopf algebras $k[Q] \to k[G]$ is injective, or equivalently if it is faithfully flat \cite[\S 14]{waterhouse}. We  recall a well-known fact.

\begin{lemma}[{\cite[Corollary 2.7]{deligne-milne}}]\label{lem:affpro}
	If $G$ is an affine group scheme over a field $k$, there exists a projective system $(G_{i})_{i}$ of group schemes of finite type over $k$ and compatible quotient homomorphisms $G \to G_{i}$ such that $G \to \projlim_{i}G_{i}$ is an isomorphism. \qed
\end{lemma}

\begin{lemma}\label{lem:nquot}
	Let $G$ be a connected, abelian affine group scheme over a field $k$ of characteristic $0$, and $n$ a positive integer prime with $\cha k$. The multiplication by $n$ homomorphism $n: G \to G$ is a quotient.
\end{lemma}

\begin{proof}
	By Lemma~\ref{lem:affpro}, $G$ is a projective limit of connected, reduced quotients $G \to Q$, and hence $k[G]$ is the union of their Hopf algebras $k[Q]$. For every such $Q$, the schematic image $nQ \subseteq Q$ is a closed subgroup with $\dim nQ = \dim Q$ thanks to \cite[Exposé XV, Lemma 1.3]{sga3.2}, which implies $nQ = Q$ since $Q$ is reduced, of finite type and irreducible \cite[\S 6.6]{waterhouse}. This implies that $n: Q \to Q$ is a quotient by \cite[\S 15.1]{waterhouse}, i.e. $n^{\#}: k[Q] \to k[Q]$ is injective. Since $k[G]$ is the union of $k[Q]$ for varying $Q$, we get that $n^{\#}: k[G] \to k[G]$ is injective as well, i.e. $n: G \to G$ is a quotient.
\end{proof}

\begin{proposition}\label{prop:gscgroup}
	Let $G$ be a connected, abelian affine group scheme over a field $k$ of characteristic $0$. 
	
	The following are equivalent.
	\begin{enumerate}
		\item $G$ is torsion free.
		\item For every positive integer $n$, the multiplication by $n$ map $n: G \to G$ is an isomorphism.
		\item For every quotient homomorphism $G \to Q$ with $Q$ of finite type and every positive integer $n$, there exists a section $G \to Q$ of the multiplication by $n$ homomorphism $n: Q \to Q$.
		\item For every quotient homomorphism $H\to Q$ with finite kernel and every homomorphism $G\to Q$, there exists a unique section $G \to H$.
	\end{enumerate}
\end{proposition}

\begin{proof}
	$(1) \Rightarrow (2)$. Since $G$ is torsion free, the multiplication by $n$ map $n: G \to G$ is a closed embedding \cite[\S 15.3]{waterhouse} and a quotient by Lemma~\ref{lem:nquot}. It follows that $n^{\#}: k[G] \to k[G]$ is bijective.
	
	$(2)\Rightarrow (1)$. Obvious.
		
	$(2)\Rightarrow (3)$. Obvious. 
	
	$(3)\Rightarrow (4)$. Write $f$ for the homomorphism $H \to Q$. Up to replacing $Q$ with the schematic image of $G \to Q$, we may assume that $G \to Q$ is a quotient and $Q$ is affine, connected, abelian. Since $f$ is affine, $H$ is affine as well. Thanks to Lemma~\ref{lem:affpro}, we may then reduce to the case in which $H,Q$ are of finite type. By replacing $H$ with the connected component of the identity, we may assume $H$ connected.
	
	After these reductions, $H$ is abelian: the commutator subgroup $[H,H]$ is connected and is a subgroup of $\ker f$, which is finite, hence $[H,H] = 0$. Since $f$ is a quotient with finite kernel and $H$ is abelian, there exists a $n$ such that $n: H \to H$ factorizes through $f$, i.e. $n = g \circ f$ for some $g: Q \to H$. The composition $f \circ g$ is the multiplication by $n$ as well: in fact, $f \circ g \circ f = f \circ n = n \circ f$, and $f$ is an epimorphism. This gives us the desired section $G \to Q \to H \to Q$.
	
	For uniqueness, we have that $H \to Q$ is étale because $\cha k = 0$ and $H \to Q$ has finite kernel. Then the diagonal is both open and closed, hence the two sections $G \to H$ coincide since $G$ is connected.
	
	$(4)\Rightarrow(2)$. Write $G=\projlim_iQ_i$ using its connected finite-type quotients. For each $i$, the multiplication map $[n]:Q_i\to Q_i$ is a quotient with finite kernel. Apply (4) to obtain the unique lift $s_i:G\to Q_i$ of the quotient map $q_i:G\to Q_i$ through $[n]$. Uniqueness makes these lifts compatible. Their limit is a homomorphism $s:G\to G$ satisfying $[n]s=\id_G$. Since $s$ commutes with multiplication by $n$, also $s[n]=\id_G$, proving (2).
\end{proof}

\begin{corollary}\label{cor:gscoh}
	Let $G$ be a connected, abelian, torsion free affine group scheme over a field $k$ of characteristic $0$. For every $i\ge 0$, $\H^{i}(k,G)$ is a $\QQ$-vector space, i.e. it is divisible and torsion free.
\end{corollary}

\begin{proof}
	Follows immediately from the fact that, for every $n$, the multiplication by $n$ map $G\to G$ is an isomorphism.
\end{proof}

A \emph{pro-torus} is a group scheme which is a projective limit of tori.

\begin{corollary}\label{cor:gsctriv}
	Let $T$ be a torsion-free pro-torus over a field $k$ of characteristic $0$, and $T \to H$ a homomorphism with $H$ of finite type. Then $\H^{1}(k,T)\to\H^{1}(k,H)$ is trivial.
	
	More precisely, there exists a quotient $T \to Q$ of finite type and a factorization $T \to Q \to H$ such that $\H^{1}(k, Q)\to\H^{1}(k, H)$ is trivial.
\end{corollary}

\begin{proof}
	Since $T$ is a pro-torus, by passing to the image we may reduce to the case in which $H$ is a torus and a quotient of $T$. 
	
	If $H$ is a torus, there exists an integer $n$ such that $\H^{1}(k, H)$ is $n$-torsion. In fact, let $k'$ be a splitting field of $H$ finite over $k$, for every $\alpha \in \H^{1}(k, H)$ we have $0 = \operatorname{cor}_{k'/k}(\alpha_{k'}) = [k':k]\alpha$, since $\alpha_{k'} = 0\in\H^{1}(k, \GG_{m}^{r})$. 
	
	Now take $Q = H$ and let $Q = H \to H$ be the multiplication by $n$ homomorphism. Clearly, $\H^{1}(k, Q) \to \H^{1}(k, H)$ is trivial. Proposition~\ref{prop:gscgroup} guarantees the existence of the factorization $T \to Q \to H$.
\end{proof}

\section{The connected component of the toric fundamental group}

We are going to give geometric conditions which guarantee that the connected component of the toric fundamental group is torsion free. Then, we will prove that this guarantees injectivity of $\cS^{\tor}_{X/k} \to \cS^{\et}_{X/k}$.

\begin{definition}
	A scheme $X$ over $k$ has \emph{virtually divisible line bundles} if, for every finite étale cover $Y\to X$, every line bundle $L$ on $Y$ and every positive integer $n$, there exists a finite étale cover $Z\to Y$ and a line bundle $M$ on $Z$ such that $L|_{Z}\simeq M^{\otimes n}$.
\end{definition}

\begin{lemma}\label{lem:curvediv}
	A smooth, projective curve of genus $\ge 1$ over a field $k$ of characteristic $0$ has virtually divisible line bundles.
\end{lemma}

\begin{proof}
	Since $\cha k = 0$, by passing to the algebraic closure, we can assume that $k$ is algebraically closed. Over an algebraically closed field, a line bundle $L$ on a smooth projective curve has the form $M^{\otimes n}$ if and only if $n\mid \deg L$. Up to passing to a finite étale cover, we can always assume that this is true if the genus is not $0$.
\end{proof}

\begin{lemma}\label{lem:abdiv}
	A torsor for an abelian variety over a field $k$ of characteristic $0$ has virtually divisible line bundles.
\end{lemma}

\begin{proof}
	Since $\cha k = 0$, then $\bar{k}/k$ is separable, hence we may easily reduce to the case in which $k = \bar{k}$. In particular, every torsor over $k$ is trivial. Let $A$ be an abelian variety and $L$ a line bundle on $A$. Every finite étale cover of $A$ is dominated, componentwise, by a multiplication map $[n]:A\to A$ for some $n>0$. Since the characteristic is $0$, for every $n$ the multiplication by $2n$ morphism $[2n]: A \to A$ is étale. Since 
	\[(2n)^{*}L \simeq L^{\otimes \frac{2n(2n+1)}{2}} \otimes (-1)^{*}L^{\otimes \frac{2n(2n-1)}{2}} =(L^{\otimes 2n+1} \otimes (-1)^{*}L^{\otimes 2n-1})^{\otimes n}\]
	by \cite[\S 6, Corollary 3]{mumford}, then $A$ has virtually divisible line bundles.
\end{proof}

\begin{lemma}\label{lem:virdiv}
	Let $X$ be a quasi-compact, quasi-separated, geometrically connected scheme over a field $k$ of characteristic $0$, and $\tilde{X} \to X$ its universal profinite étale cover. The following are equivalent.
	\begin{itemize}
		\item $X$ has virtually divisible line bundles.
		\item The group $\pic(\tilde{X})$ is divisible.
	\end{itemize}
\end{lemma}

\begin{proof}
	Notice that $\tilde{X} \to X$ is affine, hence $\tilde{X}$ is quasi-compact and quasi-separated, too. Write $\tilde{X}=\projlim_{i}X_{i}$ as a projective limit of finite étale covers. The statement follows directly from the isomorphism $\pic(\tilde{X})\simeq\indlim_{i}\pic(X_{i})$ \cite[Théorème 8.5.2]{ega-iv.3}.
\end{proof}

\begin{lemma}\label{lem:vtet}
	Let $\Upsilon$ be a vt gerbe of finite type over a field $k$. There exists a finite étale gerbe $\Upsilon^{\et}$ with a morphism 
	\[\Upsilon \to \Upsilon^{\et}\]
	which is a relative abelian gerbe whose band is a connected multiplicative group.
	
	Furthermore, if $\Upsilon \to \Upsilon^{'\et}$ is another such morphism, there exists an equivalence $\Upsilon^{\et} \to \Upsilon^{'\et}$ making the obvious diagram $2$-commutative.
	
	If $k'/k$ is an arbitrary field extension, then $(\Upsilon_{k'})^{\et} = (\Upsilon^{\et})_{k'}$.
\end{lemma}

\begin{proof}
	Define $\Upsilon^{\et} \eqdef \pet_{\Upsilon/k}$ as the étale fundamental gerbe of $\Upsilon$. If $\bar{k}$ is an algebraic closure, then $\Upsilon_{\bar{k}} = \cB G$ for some vt group $G$. By Theorem~\ref{thm:et-bc}, we get that $\Upsilon_{\bar{k}}^{\et} = \pet_{\cB G/\bar{k}} = \cB G^{\et}$, where $G^{\et}$ is the maximal étale quotient of $G$, i.e. the group scheme of the connected components. In particular, this implies that $\Upsilon \to \Upsilon^{\et}$ is a relative abelian gerbe whose band is a connected multiplicative group.
	
	If $\Upsilon \to \Upsilon^{'\et}$ is another such morphism, the universal property of the étale fundamental gerbe gives a factorization $\Upsilon \to \Upsilon^{\et} \to \Upsilon^{'\et}$. By passing to $\bar{k}$, these become morphisms $\cB G \to \cB G^{\et} \to \cB H$ for some finite group $H$; the fact that $\cB G \to \cB H$ is a relative abelian gerbe whose kernel is connected implies that $G^{\et} \to H$ is an isomorphism, and hence $\Upsilon^{\et} \to \Upsilon^{'\et}$ is an equivalence.
	
	The fact that $(\Upsilon_{k'})^{\et} = (\Upsilon^{\et})_{k'}$ follows from the uniqueness we have just proved, since $\Upsilon_{k'} \to (\Upsilon^{\et})_{k'}$ satisfies the defining property for $(\Upsilon_{k'})^{\et}$.
\end{proof}

\begin{lemma}\label{lem:pichom}
	Let $X$ be a firm scheme over a field $k$. Assume that $\pi_{1}(X_{\bar{k}})$ is trivial. The toric fundamental gerbe $\ptor_{X/k}$ is abelian banded by a connected, pro-multiplicative group scheme $T$. Furthermore,
	\[\upic_{X}(k) \simeq \hom(T,\GG_{m}).\]
	Varying $k$, this is equivalent to saying that the band $T$ is the Cartier dual $\upic_{X}^{\vee}$ of the Picard sheaf $\upic_{X}$ on the small étale site $k_{\et}$.
\end{lemma}

\begin{proof}
	Since $\pet_{X_{\bar{k}}/\bar{k}} = \spec \bar{k}$, then $\pet_{X/k} = \spec k$ by Theorem~\ref{thm:et-bc}. The fact that $\ptor_{X/k}$ is abelian banded by a connected, pro-multiplicative group scheme $T$ now follows from Lemma~\ref{lem:vtet}. In particular, $\ptor_{X/k}$ coincides with the fundamental gerbe for the class of groups of multiplicative type cf. \cite{borne-vistoli}. The isomorphism $T \simeq \upic(X)^{\vee}$ then follows from \cite[Theorem 13.11]{borne-vistoli}. Let us check this explicitely when $X(k) \neq \emptyset$ and $\pic_{X}(k) = \pic(X)$, the result then also follows by descent.
	
	A line bundle $L$ corresponds to a $\GG_{m}$-torsor, which in turn corresponds to a map $X \to \cB \GG_{m}$. By the universal property of the fundamental gerbe, this map factorizes as $X \to \ptor_{X/k} \to \cB \GG_{m}$ and hence induces a homomorphism of bands $T \to \GG_{m}$. This defines a map 
	\[\pic(X) \to \hom(T,\GG_{m}).\]
	This map is a homomorphism: in fact, if $L$, $L'$ are line bundles associated with morphisms $f,f':X \to \cB \GG_{m}$ and homomorphisms $\alpha, \alpha' \in \hom(T,\GG_{m})$, then $L \otimes L'$ is associated with the composition $X \to \cB \GG_{m} \times \cB \GG_{m} = \cB (\GG_{m} \times \GG_{m}) \to \cB \GG_{m}$, which implies that $L \otimes L'$ corresponds to $\alpha \cdot \alpha'$. Let us check injectivity and surjectivity.
	
	If the morphism of bands $T \to \GG_{m}$ induced by a line bundle $L$ is trivial, then $\ptor_{X/k} \to \cB \GG_{m}$ factorizes through $\spec k$ by Proposition~\ref{prop:canonical}, hence the same is true for $X \to \cB \GG_{m}$ and $L$ is trivial. This proves the injectivity of $\pic(X) \to \hom(T,\GG_{m})$.
	
	Let $T \to \GG_{m}$ be a homomorphism, the induced gerbe $\ptor_{X/k}\times^{T}\GG_{m}$ is banded by $\GG_{m}$ and hence is isomorphic to $\cB \GG_{m}$ since $X(k) \neq \emptyset$ implies $\ptor_{X/k}\times^{T}\GG_{m} (k) \neq \emptyset$. We thus get a morphism $X \to \cB \GG_{m}$ corresponding to a line bundle $L$. This proves surjectivity of $\pic(X) \to \hom(T,\GG_{m})$, hence we conclude.
\end{proof}

\begin{proposition}\label{prop:curvegsc}
	Let $X$ be a firm scheme over a field $k$ of characteristic $0$. The natural map
	\[\ptor_{X/k} \to \pet_{X/k}\]
	is a relative abelian gerbe whose band is a relative connected, pro-multiplicative group over $\pet_{X/k}$. In other words, the natural map from the toric fundamental group to the étale one is faithfully flat, and the kernel $T$ is a connected, pro-multiplicative group scheme.
	
	The relative band $\cT \to \pet_{X/k}$ is torsion free if and only if $X$ has virtually divisible line bundles. In other words, the  group $T$ is torsion free if and only if $X$ has virtually divisible line bundles.
\end{proposition}

\begin{proof}
	The first part follows directly from Lemma~\ref{lem:vtet}.
	
	For the second part, thanks to Theorem~\ref{thm:basechange} we may assume that $k$ is algebraically closed as well. Let $\tilde{X} \to X$ be the universal cover, it is firm by Lemma~\ref{lem:firmet}. By Proposition~\ref{prop:cartesian}, we have a $2$-cartesian diagram
	\[\begin{tikzcd}
		\tilde{X} \rar\dar	&	\ptor_{\tilde{X}/k} \rar\dar		&	\pet_{\tilde{X}/k} = \spec k \dar	\\
		X \rar				&	\ptor_{X/k} \rar					&	\pet_{X/k}.
	\end{tikzcd}\]
	
	By Lemma~\ref{lem:pichom}, $\ptor_{\tilde{X}/k}$ is abelian banded by a connected pro-multiplicative group $T$, and $\pic(\tilde{X}) = \hom(T, \GG_{m})$. The band $\cT$ of $\ptor_{X/k} \to \pet_{X/k}$ is a twisted form of $T$ over $\pet_{X/k}$, hence the relative group $\cT$ is torsion free if and only if $T$ is torsion free. Thanks to Proposition~\ref{prop:gscgroup}, this is equivalent to the following property: for every quotient homomorphism $q: T \to Q$ with $Q$ of finite type and every positive integer $n$, there exists a section $s : T \to Q$ with $n \circ s = q$.
	
	Since $T$ is pro-multiplicative, connected and $k$ is algebraically closed of characteristic $0$, then $Q = \GG_{m}^{r}$ for some $r$, hence $T$ is torsion free if and only if $\hom(T,\GG_{m}) = \pic(\tilde{X})$ is divisible. By Lemma~\ref{lem:virdiv}, this is equivalent to $X$ having virtually divisible line bundles.
\end{proof}

\begin{remark}
	With a proof analogous to the one of Lemma~\ref{lem:pichom}, one can prove that the relative band of $\ptor_{X/k} \to \pet_{X/k}$ is the Cartier dual of the relative Picard sheaf of the morphism $X \to \pet_{X/k}$.
\end{remark}

\begin{proposition}\label{prop:gsctriv2}
	Consider 
	\[f : \Upsilon \to \Phi\]
	a locally full morphism of pro-algebraic gerbes over a field $k$ of characteristic $0$.
	
	Assume that $f$ is a relative abelian gerbe whose band is a relative torsion free pro-torus over $\Phi$. In other words, for any geometric section of $\Upsilon$ assume that the induced homomorphism of automorphism group schemes is faithfully flat and its kernel is a torsion free pro-torus.
	
	The map 
	\[\Upsilon(k)_{\rmq} \to \Phi(k)_{\rmq}\]
	is injective.
	
	Furthermore, if $s \in \Phi(k)_{\rmq}$ is not in the image, there exists a $2$-commutative diagram of gerbes
	\[\begin{tikzcd}
		\Upsilon \rar\dar	&	\Phi\dar	\\
		\Upsilon'\rar		&	\Phi'
	\end{tikzcd}\]
	such that every morphism is locally full, $\Upsilon'$ and $\Phi'$ are of finite type over $k$, $\Upsilon' \to \Phi'$ is a relative gerbe banded by a torus and the image of $s$ in $\Phi'$ does not lift to $\Upsilon'$.
\end{proposition}

\begin{proof}
	Write $\Upsilon = \projlim_{i} \Upsilon_{i}$ as a projective limit of gerbes of finite type over $k$, with locally full projections. For every $i$, consider the inertia stack $\rI_{\Upsilon_{i}} \to \Upsilon_{i}$. Recall that inertia stacks are relative groups \cite[\href{https://stacks.math.columbia.edu/tag/050R}{Tag 050R}]{stacks-project}. There is a natural $2$-commutative diagram
	\[\begin{tikzcd}
		\rI_{\Upsilon} \rar \dar	&	\rI_{\Upsilon_{i}} \dar		\\
		\Upsilon \rar				&	\Upsilon_{i}
	\end{tikzcd}\]
	where the upper horizontal arrow is a homomorphism of relative groups. The inertia $\rI_{f}$ of $f$ defines a subgroup of $\rI_{\Upsilon}$. Define $\Phi_{i}$ as the rigidification \cite[Appendix A]{abramovich-olsson-vistoli} of $\Upsilon_{i}$ by the image of $\rI_{f}$ in $\rI_{\Upsilon_{i}}$.
	
	Since $f$ is locally full and $\rI_{f}$ maps to $0$ in the relative inertia of $\Upsilon \to \Phi_{i}$, there is a natural factorization $\Upsilon \to \Phi \to \Phi_{i}$ (the proof is analogous to the one of \cite[Theorem 5.1.5(2)]{abramovich-corti-vistoli}).	The natural morphism $\Phi \to \projlim_{i}\Phi_{i}$ is an isomorphism, and $\Upsilon_{i} \to \Phi_{i}$ is a relative abelian gerbe whose band is a relative torus. In fact, this can be checked after base changing to $\bar{k}$, where it reduces to the following elementary fact: if $G = \projlim_{i} G_{i}$ is a projective limit of algebraic groups $G_{i}$ of finite type over $\bar{k}$, $G \to H$ is a quotient homomorphism whose kernel $T$ is a pro-torus, and $T_{i}$ is the image of $T$ in $G_{i}$, then $T_{i}$ is a torus and $H = \projlim_{i} G_{i}/T_{i}$.
	
	Let $(\varphi_{i})_{i}\in\Phi(k)_{\rmq}=\lim_{i}\Phi_{i}(k)/\unsim$ be a quasi-section of $\Phi$. It is enough to prove the following: if $\varphi_{i}$ lifts to $\Upsilon_{i}$ for every $i$, then $(\varphi_{i})_{i}$ lifts uniquely to $\Upsilon(k)_{\rmq}$. 
	
	Let $\Psi_{i}\to\Upsilon_{i}$ be the inverse image of $\varphi_{i}:\spec k \to \Phi_{i}$, then $\Psi_{i}$ is an abelian gerbe banded by a torus $T_{i}$. For every $j \ge i$, a choice of an isomorphism between the image of $\varphi_{j}$ and $\varphi_{i}$ in $\Phi_{i}(k)$ gives a morphism $\Psi_{j}\to\Psi_{i}$. These morphisms, including the induced maps on bands, may depend on the chosen isomorphisms; we do not need compatible choices. 
	
	For every $i$, $\Psi_{i}(k)/\unsim$ is a principal homogeneous space for $\H^{1}(k, T_{i})$. The argument of Corollary~\ref{cor:gsctriv}, applied to the relative torsion free pro-torus band of $\Upsilon\to\Phi$ and using finite presentation, shows that, for $j>>i$ large enough and any chosen comparison, $\H^{1}(k,T_j)\to\H^{1}(k,T_i)$ is trivial, and hence the image of $\Psi_{j}(k)/\unsim \to \Psi_{i}(k)/\unsim$ contains exactly one element. While this element may depend on the choice of the isomorphism, its image $\upsilon_{i}$ in $\Upsilon_{i}(k)/\unsim$ does not, since the map $\Psi_j(k)/\unsim\to\Upsilon_i(k)/\unsim$ does not depend on the choice. Comparing two sufficiently large indices with a common upper index shows that $\upsilon_i$ does not depend on $j >> i$.
	
	Clearly, if $j > i$ then $\upsilon_{j}$ maps to $\upsilon_{i}$, and $(\upsilon_{i})_{i}$ defines an element of $\Upsilon(k)_{\rmq}$ which lifts $(\varphi_{i})_{i}$.
	
	If $(\upsilon_{i}')_{i}$ is another element mapping to $(\varphi_{i})_{i}$, then by definition of $\Psi_{i}$ we can choose a further lifting $\psi_{i}' \in \Psi_{i}(k)/\unsim$ of $\upsilon_{i}'$ for every $i$. For $j >> i$ large enough, the image of $\psi_j'$ in $\Upsilon_i(k)/\unsim$ is $\upsilon_i$ by its defining property. This implies that $\upsilon_{j}'$ maps to $\upsilon_{i}$, and hence $\upsilon_{i}' = \upsilon_{i}$. This concludes the proof.
\end{proof}

\begin{theorem}\label{thm:toretinj}
	Let $X$ be a firm scheme over a field $k$ of characteristic $0$. If $X$ has virtually divisible line bundles, the natural map $\cS^{\tor}_{X/k} \to \cS^{\et}_{X/k}$ is injective.
\end{theorem}

\begin{proof}
	This is a direct consequence of Propositions~\ref{prop:curvegsc} and \ref{prop:gsctriv2}.
\end{proof}

\section{Toric Galois sections and points of the Picard scheme}

Let $X$ be a proper, geometrically connected variety over a field $k$ of arbitrary characteristic. If $X(k) \neq \emptyset$, then $\pic(X) = \upic_{X}(k)$. In general, we only have an inclusion $\pic(X) \subseteq \upic_{X}(k)$, and there is a cohomological obstruction $\upic_{X}(k) \to \br(k)$ whose kernel is the Picard group $\pic(X)$ \cite[Remark 9.2.11]{kleiman}. The main example is when $X$ is a non-trivial Brauer--Severi variety, so that $\cO(1) \in \upic_{X}(k)$ is not a line bundle. For such a proper variety $X$, the relative Brauer group $\br(X/k)$ coincides with the image of the cohomological obstruction, hence $\br(X/k) = \upic_{X}(k)/\pic(X)$. 

The Leray spectral sequence in fppf cohomology for $X\to\spec k$ gives the exact sequence
\[
0\to\pic(X)\to\upic_X(k)\to\br(k)\to\H^2_{\mathrm{fppf}}(X,\GG_m),
\]
see \cite[Remark 9.2.11 and Exercise 9.3.11]{kleiman}.

In particular, if $\pic(X) \neq \upic_{X}(k)$ then $\br(X/k) \neq 0$, and $\cS^{\tor}_{X/k} = \emptyset$ by Lemma~\ref{lem:brauerobstruction}. 

On the other hand, J. Stix has proved that, if $p \in \upic_{X}(k)$ is a torsion point of order prime with $\cha k$ and $p \not \in \pic(X)$, then $\cS^{\et}_{X/k} = \emptyset$ \cite[Proposition 12]{stix-period}. 

The following Proposition~\ref{prop:obstruction} gives a unified, alternative proof of these two statements. Even if it's not strictly necessary, we find it illuminating because it replaces the long exact sequence above with a geometric argument based on the Picard stack.

\begin{proposition}\label{prop:obstruction}
	Let $X$ be a proper, geometrically connected variety over a field $k$. 
	
	A rational point $p \in \upic_{X}(k)$ defines a map $X \to \cG$ such that
	\begin{itemize}
		\item $\cG$ is a gerbe banded by $\GG_{m}$,
		\item $\cG(k) \neq \emptyset$ if and only if $p \in \pic(X) \subseteq \upic_{X}(k)$,
		\item if $k'$ is any field extension with $X(k') \neq \emptyset$, the induced morphism $X_{k'} \to \cG_{k'} = B_{k'} \GG_{m}$ defines the $\GG_{m}$-torsor associated with the line bundle $p_{k'} \in \pic(X_{k'}) = \upic_{X_{k'}}(k')$,
		\item if $p$ is torsion of order prime with $\cha k$, we have a factorization $X \to \Pi_{X/k}^{\et} \to \cG$.
	\end{itemize}
		
	As a consequence, if $\pic(X) \neq \upic_{X}(k)$ then $\cS^{\tor}_{X/k} = \emptyset$, and if there is a torsion element $p \in \upic_{X}(k) \setminus \pic(X)$ of order prime with $\cha k$, then $\cS^{\et}_{X/k} = \emptyset$.
\end{proposition}

\begin{proof}
	Recall that the Picard stack $\cP ic_{X}$ is defined as follows: if $S$ is a scheme over $k$, then $\cP ic_{X}(S)$ is the groupoid of line bundles on $X_{S}$, so that $\pic(X_{S})$ is the set of isomorphism classes of $\cP ic_{X}(S)$ \cite[\href{https://stacks.math.columbia.edu/tag/0372}{Tag 0372}]{stacks-project}. The Picard sheaf $\upic_{X}$ is the fppf sheafification of this stack; under our hypotheses, it is representable by a group scheme \cite[Corollary 9.4.18.3]{kleiman}. 
	
	A Picard functor is usually defined using the quotients $\pic(X_{S}) / \pic(S)$, and then sheafified: notice however that $\pic(S)$ is automatically killed by the sheafification, since every line bundle over $S$ is locally trivial in $S$. The reason for considering the functor $S \mapsto \pic(X_{S}) / \pic(S)$ is that, if $X(k) \neq \emptyset$, the functor is already an fppf sheaf \cite[Theorem 9.2.5]{kleiman}, thus simplifying things considerably.
	
	There is a natural morphism $\cP ic_{X} \to \upic_{X}$ which is a gerbe banded by $\GG_{m}$. In fact, given any scheme $S$ and a line bundle $L$ on $X_{S}$, the group of automorphisms of $L$ is naturally $\GG_{m}(S) = \cO(S)^{*}$ because $\cO(X_{S}) = \cO(S)$ (recall that $X$ is proper, geometrically connected). If $p \in \upic_{X/k}(k)$ is a rational point, we write $\cG_{p} \subset \cP ic_{X}$ for the fiber, it is a gerbe over $k$ banded by $\GG_{m}$. The gerbe $\cG_{p}$ has a rational point if and only if $p \in \upic_{X}(k)$ defines a line bundle.
	
	Write $S = X$ for a copy of $X$. The fibered product $X_{S} = X \times X$ has a section $s: S \to X_{S}$, the diagonal. By \cite[Theorem 9.2.5]{kleiman}, this implies that there is a line bundle $M$ on $X_{S}$ corresponding to $p \in \upic_{X}(k)$. The line bundle $M$ is only well defined in $\pic(X_{S})/\pic(S)$ (because there are several morphisms $S \to \cG_{p}$); since we can modify $M$ by elements of $\pic(S)$, we can rigidify the situation by imposing that $M$ restricts to the trivial bundle on the section $s: S \to X_{S}$.
	
	The line bundle $M$ on $X_{S}$ corresponds to a morphism $S \to \cG_{p}$; this factorizes through $\ptor_{S/k}$, because $\cG_{p}$ is banded by $\GG_{m}$. Since $S = X$ and $\cG_{p}(k)$ is empty if $p \notin \pic(X)$, we get that $\cS^{\tor}_{X/k} = \emptyset$ if $\upic_{X}(k) \neq \pic(X)$.
	
	It remains to give the factorization $X \to \pet_{X/k} \to \cG_{p}$ if $p$ is torsion of order prime with $\cha k$. By Proposition~\ref{prop:canonical}, we have a factorization $\ptor_{S/k} \to \Phi_{p} \to \cG_{p}$, with $\Phi_{p} \to \cG_{p}$ faithful. Since $\cG_{p}$ is banded by $\GG_{m}$, either $\Phi_{p} \simeq \cG_{p}$ or $\Phi_{p}$ is banded by $\mu_{n} \subset \GG_{m}$ for some $n$. In the second case, if $\cha k \nmid n$ then we get the desired factorization by the universal property of the étale fundamental gerbe.
		
	I claim that $\Phi_{p}$ is banded by $\mu_{n}$ if and only if $p \in \upic_{X}(k)$ is torsion of order $n$. This is enough to conclude, and can be checked after extending the base field. Hence, we may assume that $X(k) \neq \emptyset$, so that applying \cite[Theorem 9.2.5]{kleiman} again we get a line bundle $L \in \pic(X) = \upic_{X}(k)$ associated with $p$ and a section $\spec k \to \cG_{p}$ (this section is unique up to isomorphism, because $\H^{1}(k,\GG_{m})$ is trivial).
	
	Hence, we have a map $X \to B\GG_{m}$ corresponding to the line bundle $L \in \pic(X) = \H^{1}(X, \GG_{m})$. By the Kummer exact sequence
	\[\H^{1}(X, \mu_{n}) \to \H^{1}(X, \GG_{m}) \xrightarrow{n} \H^{1}(X, \GG_{m}),\]
	if $L$ is $n$-torsion then it is induced by an element of $\H^{1}(X, \mu_{n})$. The corresponding map $X \to B\mu_{n}$ gives the desired factorization.
\end{proof}

\begin{lemma}\label{lem:bstor}
	Let $P$ be a Brauer--Severi variety of positive dimension over a field $k$. The gerbe $\ptor_{P/k}$ is of finite type, abelian and banded by $\GG_{m}$. If $P(k)=\emptyset$, then $\ptor_{P/k}(k)=\emptyset$.
\end{lemma}

\begin{proof}
	Let $P \to \cG$ be the morphism given by Proposition~\ref{prop:obstruction} for $\cO(1) \in \upic_{P}(k)$. The gerbe $\cG$ is abelian and banded by $\GG_{m}$, hence there is an induced morphism $\ptor_{P/k} \to \cG$. This is an isomorphism: it can be checked after passing to a separable closure so that $P = \PP^{n}$, where it follows from the fact that $\PP^{n} \to \cG = B \GG_{m}$ corresponds to $\cO(1)$, plus Lemma~\ref{lem:pichom}. Finally, $\cG(k) \neq \emptyset$ if and only if $\cO(1) \in \pic(P)$, which in turn is equivalent to $P(k) \neq \emptyset$.
\end{proof}

\begin{corollary}\label{cor:brtor}
	If $X$ is a firm fibered category over a field $k$ with a morphism $X \to P$, where $P$ is a non-trivial Brauer--Severi variety, then $\cS^{\tor}_{X/k}=\emptyset$.
\end{corollary}

\begin{proof}
	We have an induced map $\cS^{\tor}_{X/k}\to\cS^{\tor}_{P/k}$; since $\ptor_{P/k}$ is of finite type, it is algebraic, hence $\cS_{P/k}^{\tor}=\ptor_{P/k}(k)_{\rmq}=\ptor_{P/k}(k)/\unsim=\emptyset$ by Lemma~\ref{lem:bstor}.
\end{proof}	

\section{Toric versus étale}

In this section, we compare the classical version of the section conjecture with the toric one.

\begin{lemma}\label{lem:torsur}
	Let $X$ be a firm fibered category over a field $k$ and $Y \to X$ a geometrically connected profinite étale cover. The natural map
	\[\cS^{\tor}_{Y/k} \to \cS^{\et}_{Y/k} \times_{\cS^{\et}_{X/k}} \cS^{\tor}_{X/k}\]
	is surjective.
\end{lemma}

\begin{proof}
	By Proposition~\ref{prop:cartesian}, we have a $2$-cartesian diagram
	\[\begin{tikzcd}
		\ptor_{Y/k} \rar \dar 	&	\pet_{Y/k} \dar	\\
		\ptor_{X/k} \rar		&	\pet_{X/k}
	\end{tikzcd}\]
	which we may write as a projective limit of $2$-cartesian diagrams
	\[\begin{tikzcd}
		\Pi \rar \dar 		&	\Pi^{\et} \dar		\\
		\Upsilon \rar		&	\Upsilon^{\et}
	\end{tikzcd}\]
	with $\Pi$, $\Upsilon$ vt gerbes of finite type. The map
	\[\Pi(k)/\unsim \to \Upsilon(k)/\unsim \times_{\Upsilon^{\et}(k)/\unsim} \Pi^{\et}(k)/\unsim\]
	is easily checked to be surjective. Furthermore, its fibers are finite, since two elements of $\Upsilon^{\et}(k)$ can only be isomorphic in a finite number of ways.
	
	The limit of the left hand side is $\cS^{\tor}_{Y/k}$, and the limit of the right hand side is $\cS^{\et}_{Y/k} \times_{\cS^{\et}_{X/k}} \cS^{\tor}_{X/k}$. A projective limit of surjective maps with finite fibers is surjective, hence we conclude.
\end{proof}

\begin{corollary}\label{cor:torbij}
	Let $X$ be a firm scheme with virtually divisible line bundles over a field $k$ of characteristic $0$, and $Y \to X$ a geometrically connected profinite étale cover. The natural map
	\[\cS^{\tor}_{Y/k} \to \cS^{\et}_{Y/k} \times_{\cS^{\et}_{X/k}} \cS^{\tor}_{X/k}\]
	is bijective, i.e. $\cS^{\tor}_{Y/k} \subset \cS^{\et}_{Y/k}$ is the inverse image of $\cS^{\tor}_{X/k} \subset \cS^{\et}_{X/k}$. \qed
\end{corollary}

Given a scheme $X$ over a field $k$ and $T$ a torus, write $\H^{2}(X/k, T)$ for the kernel of $\H^{2}(k, T) \to \H^{2}(X, T)$.

Recall that an étale neighbourhood of a section $s \in \cS^{\et}_{X/k}$ is a geometrically connected finite étale cover $Y \to X$ together with a section $r \in \cS^{\et}_{Y/k}$ whose image is $s$.

\begin{lemma}\label{lem:torusdiv}
	Let $X$ be a firm scheme with virtually divisible line bundles over a field $k$ of characteristic $0$ with a Galois section $s \in \cS^{\et}_{X/k}$, $T$ a torus over $k$ and $\psi \in \H^{2}(X/k, T)$. 
	
	For every positive integer $d$, there exists an étale neighbourhood $Y \to X$ of $s$ and a cohomology class $\psi/d \in \H^{2}(Y/k, T)$ such that $d \cdot (\psi/d) = \psi$.
\end{lemma}

\begin{proof}
	Write $\Psi$ for the gerbe banded by $T$ corresponding to $\psi$, the fact that $\psi$ maps to $0 \in \H^{2}(X, T)$ is equivalent to the existence of a morphism $X \to \Psi$.
	
	Denote by $X_{s}$ the decomposition tower of $s$; that is, the projective limit of the étale neighbourhoods of $s$, which coincides with the fibered product $X \times_{\pet_{X/k}} \spec k$ by Proposition~\ref{prop:cartesian}. The toric fundamental gerbe $\ptor_{X_{s}/k}$ of $X_{s}$ is an abelian gerbe banded by a torsion free pro-torus by Proposition~\ref{prop:curvegsc}. Write $\Theta$ for the band, we have an induced morphism $\ptor_{X_{s}/k} \to \Psi$ and hence a homomorphism $f : \Theta \to T$. Since $\Theta$ is torsion free, by Proposition~\ref{prop:gscgroup} we have a factorization $\Theta \xrightarrow{f/d} T \xrightarrow{d} T$ of $f$.
	
	Write $\Psi/d$ for the $T$-gerbe induced by $\ptor_{X_{s}/k}$ along $f/d$, we have induced morphisms $X_{s} \to \Psi/d \to \Psi$. If $\psi/d \in \H^{2}(k,T)$ is the cohomology class associated with $\Psi/d$, then we have $d \cdot (\psi/d) = \psi$ since the morphism $\Psi/d \to \Psi$ corresponds to the homomorphism of bands $d : T \to T$. Since $\Psi/d$ is of finite type and $X_{s}$ is the limit of the étale neighbourhoods of $s$, there exists an étale neighbourhood $Y \to X$ with a morphism $Y \to \Psi/d$, hence $\psi/d \in \H^{2}(Y/k, T)$.
\end{proof}

\begin{corollary}\label{cor:torusbrauer}
	Let $X$ be a firm scheme with virtually divisible line bundles over a field $k$ of characteristic $0$, $s \in \cS^{\et}_{X/k}$ a Galois section, and $T$ a torus.
	
	If $\H^{2}(X/k, T) \neq 0$, there exists an étale neighbourhood $Y \to X$ of $s$ and a finite extension $k'/k$ with $\br(Y_{k'}/k') \neq 0$.
\end{corollary}

\begin{proof}
	Let $\psi \in \H^{2}(X/k, T)$ be a non-trivial element. Fix $k'/k$ a splitting field for $T$ finite over $k$, that is, $T_{k'} \simeq \GG_{m}^{r}$. Write $d$ for the degree of $k'/k$. 
	
	By Lemma~\ref{lem:torusdiv}, there exists an étale neighbourhood $Y \to X$ of $s$ and a cohomology class $\psi/d \in \H^{2}(Y/k, T)$ such that $d \cdot (\psi/d) = \psi$.
	
	Notice that $(\psi/d)_{k'} \in \H^{2}(Y_{k'}/k',T_{k'}) = \br(Y_{k'}/k')^{r}$. We have 
	\[\operatorname{cor}_{k'/k}\left((\psi/d)_{k'} \right) = d \cdot (\psi/d) = \psi \neq 0,\]
	hence $(\psi/d)_{k'} \neq 0$. This concludes the proof.
\end{proof}

\begin{proposition}\label{prop:toricbrauer}
	Let $X$ be a firm scheme with virtually divisible line bundles over a field $k$ of characteristic $0$.
	
	The subset $\cS^{\tor}_{X/k} \subseteq \cS^{\et}_{X/k}$ coincides with the set of étale Galois sections $s$ such that, for every étale neighbourhood $Y \to X$ and every finite extension $k'/k$, the relative Brauer group $\br(Y_{k'}/k')$ is trivial.
\end{proposition}

\begin{proof}
	Write $B \subset \cS^{\et}_{X/k}$ for the subset of Galois sections described in the statement. 
	
	Let $s \in \cS^{\tor}_{X/k}$ be a toric Galois section, and $Y \to X$, $r \in \cS^{\et}_{Y/k}$ an étale neighbourhood of the corresponding étale Galois section, we have that $r \in \cS^{\tor}_{Y/k}$ is toric thanks to Corollary~\ref{cor:torbij}. In particular, $\cS^{\tor}_{Y/k} \neq \emptyset$, hence $\br(Y_{k'}/k')$ is trivial for any finite extension $k'$ thanks to Lemma~\ref{lem:brauerobstruction}. This implies that $\cS^{\tor}_{X/k} \subset B$.
	
	On the other hand, if $s \in \cS^{\et}_{X/k}$ is an étale Galois section not in $\cS^{\tor}_{X/k}$, by Proposition~\ref{prop:gsctriv2} we have a $2$-commutative diagram
	\[\begin{tikzcd}
		\ptor_{X/k} \rar \dar	&	\pet_{X/k} \dar	\\
		\Upsilon \rar			&	\Phi
	\end{tikzcd}\]
	of locally full morphisms with $\Upsilon$, $\Phi$ of finite type over $k$, $\Upsilon \to \Phi$ a relative abelian gerbe banded by a torus, and the image $\varphi$ of $s$ in $\Phi(k)$ does not lift to $\Upsilon (k)$. Define $Y$, $\Psi$ by the $2$-cartesian diagram
	\[\begin{tikzcd}
		Y \rar \dar		&	\Psi \rar \dar	&	\spec k \dar["\varphi"]	\\
		X \rar			&	\Upsilon \rar	&	\Phi.
	\end{tikzcd}\]
	Since $\varphi$ is the image of $s$, then $Y$ is an étale neighbourhood of $s$ and $\Psi$ is an abelian gerbe banded by a torus $T$ over $k$ with $\Psi(k) = \emptyset$. The cohomology class of $\Psi$ then defines a non-trivial element of $\H^{2}(Y/k, T)$, hence $s \notin B$ thanks to Corollary~\ref{cor:torusbrauer}. This proves that $B \subset \cS^{\tor}_{X/k}$, and concludes the proof.
\end{proof}

\begin{theorem}\label{thm:toretgeneral}
	Let $k_{0}$ be a field of characteristic $0$. The following are equivalent.
	\begin{itemize}
		\item For every smooth, projective hyperbolic curve defined over a field $k$ finite over $k_{0}$, the profinite Kummer map $X(k) \to \cS^{\et}_{X/k}$ is bijective.
		\item For every smooth, projective hyperbolic curve defined over a field $k$ finite over $k_{0}$, the toric Kummer map $X(k) \to \cS^{\tor}_{X/k}$ is bijective, and if $\br(X/k) \neq 0$ then $\cS^{\et}_{X/k}=\emptyset$.
	\end{itemize}
\end{theorem}

\begin{proof}
	Assume the first condition, and let $X/k$ be as in the statement. Since $X(k) \to \cS^{\et}_{X/k}$ is bijective by assumption and $\cS^{\tor}_{X/k} \to \cS^{\et}_{X/k}$ is injective by Theorem~\ref{thm:toretinj}, then $X(k) \to \cS^{\tor}_{X/k}$ is bijective. Furthermore, if $\br(X/k) \neq 0$ then $X(k) = \emptyset$ and hence $\cS^{\et}_{X/k} = \emptyset$ by hypothesis.
	
	Assume the second condition, and let $X/k$ be as in the statement. We want to prove that $X(k) \to \cS^{\et}_{X/k}$ is bijective; since $X(k) \to \cS^{\tor}_{X/k}$ is bijective by assumption, we reduce to proving that $\cS^{\tor}_{X/k} \to \cS^{\et}_{X/k}$ is bijective. By Proposition~\ref{prop:toricbrauer}, this amounts to proving that, for every étale Galois section $s \in \cS^{\et}_{X/k}$, every étale neighbourhood $Y \to X$ and every finite extension $k'/k$, the relative Brauer group $\br(Y_{k'}/k')$ is trivial. This is guaranteed by the assumption, hence we conclude.
\end{proof}

\section{Abelian varieties over $p$-adic fields}

We want now to prove Theorem~\ref{thm:torab}, i.e. the toric section conjecture for abelian varieties over $p$-adic fields, and use this to give a different point of view on toric Galois sections of curves. More generally, we prove the following.

\begin{theorem}\label{thm:torab+}
	Let $B$ be a torsor for an abelian variety $A$ over a field $k$ finite over $\QQ_{p}$. Then
	\[B(k) \to \cS^{\tor}_{B/k}\]
	is bijective.
\end{theorem}

Let $B$ be a torsor for an abelian variety $A$ over a field $k$ finite over $\QQ_{p}$. By Lemma~\ref{lem:abdiv} and Theorem~\ref{thm:toretinj}, the map
\[\cS^{\tor}_{B/k} \to \cS^{\et}_{B/k}\]
is injective. Furthermore, 
\[B(k) \to \cS^{\tor}_{B/k}\]
is injective as well, because $B(k) \to \cS^{\et}_{B/k}$ is injective \cite[Proposition 73]{stix}.

Let 
\[\widehat{A} = \upic^{0}_{A/k} = \upic^{0}_{B/k}\]
be the dual abelian variety. Recall that $\widehat{A}(k)$ identifies with $\ext^{1}_{k}(A,\GG_{m})$ \cite[Lemma 3.1]{milne-duality}. We have the local Tate duality
\[\H^{1}(k,A) \times \widehat{A}(k) \to \H^{2}(k,\GG_{m}) = \br(k) = \QQ/\ZZ,\]
i.e. a pairing which identifies the torsion group $\H^{1}(k,A)$ with the Pontryagin dual $\hom(\widehat{A}(k), \QQ/\ZZ)$, and the profinite group $\widehat{A}(k)$ with the Pontryagin dual $\hom(\H^{1}(k,A), \QQ/\ZZ)$, see \cite[Theorem 3.2]{milne-duality}.

Furthermore, there is an exact sequence
\[0 \to \pic(B) \to \upic_{B}(k) \xrightarrow{b} \br(k) = \QQ/\ZZ \to \br(B),\]
see \cite[\S 2]{stix-period}, where $b:\upic_{B}(k) \xrightarrow{b} \br(k)$ is called the \emph{Brauer obstruction}.

\begin{lemma}\label{lem:tatebrauer}
	Let $B$ be a torsor for an abelian variety $A$ over a field $k$ finite over $\QQ_{p}$. The Brauer obstruction
	\[b:\upic_{B}(k) \to \br(k)\]
	restricted to $\upic^{0}_{B/k} = \widehat{A}$ is the local Tate pairing, i.e.
	\[b(x) = ([B],x)\]
	for $x \in \widehat{A}(k)$ (up to a sign).
\end{lemma}

\begin{proof}
	Let $k'/k$ be a finite Galois extension and write $\Gamma = \gal(k'/k)$. We have that $\widehat{A}(k) = \ext^{1}_{k}(A,\GG_{m})$. The Tate pairing
	\[\H^{1}(\Gamma, A(k')) \times \ext^{1}_{k}(A,\GG_{m}) \to \H^{2}(\Gamma, k'^{*})\]
	is defined as follows. Given a $1$-cocycle $\gamma: \Gamma \to A(k')$ and a $\Gamma$-invariant extension $E$ of $A(k')$ by $k'^{*}$, fix an automorphism $g: E \to E$ commuting with the Galois actions on $A(k')$ and $k'^{*}$ for every $g \in \Gamma$ and lift $\gamma$ to a function $\gamma': \Gamma \to E$. The Tate pairing is the cohomology class of the cocycle $\lambda: \Gamma^{2} \to k'^{*}$
	\[\lambda(g,h) = \gamma'(gh) \cdot g\gamma'(h)^{-1} \cdot \gamma'(g)^{-1} \in k'^{*} \subset E.\]
	
	Now assume that $B(k') \neq \emptyset$, and let $L$ be a Galois invariant, algebraically trivial line bundle on $B_{k'}$. The Brauer obstruction $b(L)$ is defined as follows \cite[\S 2]{stix-period}: for any $g \in \Gamma$, fix $\phi(g)$ an automorphism of $L$ over $g^{\#}:B_{k'} \to B_{k'}$, then $b(L)$ is the cohomology class of the cocycle 
	\[\omega(g,h) = \phi(gh) \circ \phi(h)^{-1} \circ \phi(g)^{-1} \in \aut_{k'}(L) = k'^{*}.\]
	
	Choose $\spec k' \to B$ a $k'$-point and identify it with the origin of $A(k') = B(k')$, its Galois conjugates define a cocycle $\gamma : \Gamma \to B(k') = A(k')$ such that $[\gamma] = [B]$. If $e:\spec k' \to L^{*}$ lifts $\spec k' \to B$ and we use it as origin of $L^{*}(k') $, the automorphism $\phi(g)$ of $L$ defines a group automorphism $g: L^{*}(k') \to L^{*}(k')$, $g(p) = \phi(g)(p) \cdot \phi(g)(e)^{-1}$ commuting with the Galois actions on $A(k')$ and $k'^{*}$.
	
	We can define a function $\gamma' : \Gamma \to E = L^{*}(k')$, $\gamma'(g) = \phi(g)(e)$ lifting $\gamma$. Since $\omega(g,h)$ coincides with its evaluation at any point, we get that
	\[\omega(g,h)(\phi(g)(\phi(h)(e))) = \phi(gh) \circ \phi(h)^{-1} \circ \phi(g)^{-1}(\phi(g)(\phi(h)(e))) = \]
	\[ = \phi(gh)(e) = \phi(gh)(e) \cdot \phi(g)(\phi(h)(e))^{-1} \cdot \phi(g)(\phi(h)(e))\]
	and hence
	\[ \omega(g,h) = \phi(gh)(e) \cdot \phi(g)(\phi(h)(e))^{-1} = \]
	\[ = \gamma'(gh) \cdot \phi(g)(\gamma'(h))^{-1} = \gamma'(gh) \cdot g\gamma'(h)^{-1} \cdot \gamma'(g)^{-1} = \lambda(g,h).\]
	This concludes the proof.
\end{proof}

\begin{corollary}\label{cor:abbr}
	The relative Brauer group $\br(B/k)$ is trivial if and only if $B(k) \neq \emptyset$.
\end{corollary}

\begin{proof}
	If $B(k) = \emptyset$, then $[B] \in \H^{1}(k,A)$ defines a non-trivial homomorphism $\widehat{A}(k) \to \QQ/\ZZ$ by local Tate duality. Let $a \in \widehat{A}(k)$ be any point such that $([B], a) \neq 0$, then $a \in \widehat{A}(k) = \upic^{0}_{B/k}(k)$ is not in $\pic(B)$ by Lemma~\ref{lem:tatebrauer}.
\end{proof}

\begin{lemma}\label{lem:abnb}
    Let $A$ be an abelian variety over a field $k$ finite over $\QQ_{p}$, and $s \in \cS^{\et}_{A/k}$ a Galois section of $A$. If $s$ is not geometric, there exists an étale neighbourhood $B\to A$ such that $B(k)=\emptyset$.
\end{lemma}

\begin{proof}
    Let $\tilde{A}=\projlim_{i}A_{i}$ be the projective limit of the étale neighbourhoods $A_{i}$ of $s$, i.e. the decomposition tower of $s$ in the terminology of \cite{stix}. Since $s$ is not geometric, then $\tilde{A}(k)=\emptyset$ \cite[Lemma 52]{stix}. Since $k$ is finite over $\QQ_{p}$, then $A_{i}(k)$ is compact for the $p$-adic topology; the fact that $\tilde{A}(k)=\projlim_{i}A_{i}(k)$ is empty then implies that $A_{i}(k)=\emptyset$ for some $i$.
\end{proof}

\begin{proof}[Proof of Theorem~\ref{thm:torab+}] 
	Let $s \in \cS^{\tor}_{B/k}$ be a toric Galois section. If, by contradiction, $s$ is not geometric, thanks to Lemma~\ref{lem:abnb} we can pass to an étale neighbourhood and assume that $B(k) = \emptyset$, so that $B$ is a non-trivial torsor for an abelian variety $A$. By Corollary~\ref{cor:abbr}, $\br(B/k)$ is non-trivial, hence $\cS^{\tor}_{B/k} = \emptyset$ by Lemma~\ref{lem:brauerobstruction}, giving a contradiction. This concludes the proof of Theorem~\ref{thm:torab+}.
\end{proof}

Let us see the consequences of Theorem~\ref{thm:torab} for hyperbolic curves. We recall a result of J. Stix \cite[Theorem 15, 16]{stix-period}.

\begin{theorem}[Stix]\label{thm:stix}
    Let $X$ be a smooth projective curve of genus $\ge 1$ over a field $k$ which is a finite extension of $\QQ_{p}$. If $X$ has a Galois section $s$, the index of $X$ is a power of $p$, and if it is different from $1$ then $s$ has an étale neighbourhood $Y \to X$ with $\upic^{1}_{Y}(k) = \emptyset$.
\end{theorem}

\begin{proposition}\label{prop:torindex}
    Let $X$ be a smooth projective curve of genus $\ge 1$ over a field $k$ which is a finite extension of $\QQ_{p}$, and $s\in\cS^{\et}_{X/k}$ a Galois section. The following are equivalent.
    \begin{enumerate}
        \item For every étale neighbourhood $(Y,r)$ of $s$, the image $r^{\rm ab}$ of $r$ in $\cS^{\et}_{\upic^{1}_{Y}/k}$ is geometric.		
        \item For every étale neighbourhood $(Y,r)$ of $s$, $\upic^{1}_{Y}(k)\neq\emptyset$.
        \item For every étale neighbourhood $(Y,r)$ of $s$, $Y$ has index $1$.
        \item $s \in \cS^{\tor}_{X/k} \subset \cS^{\et}_{X/k}$.
    \end{enumerate}
\end{proposition}

\begin{proof}
	$(1) \Rightarrow (2)$ Obvious.
	
	$(2) \Rightarrow (3)$ Follows from Stix' Theorem~\ref{thm:stix}.
	
	$(3) \Rightarrow (1)$ If, by contradiction, the image $r^{\rm ab}$ of $r$ in $\cS^{\et}_{\upic^{1}_{Y}/k}$ is not geometric, then by Lemma~\ref{lem:abnb} it has an étale neighbourhood $B \to \upic^{1}_{Y}$ with $B(k) = \emptyset$. If $Z \eqdef B \times_{\upic^{1}_{Y}} Y$, then $Z$ is an étale neighbourhood of $r$ and we have an induced map $\upic^{1}_{Z} \to B$ since $\upic^{1}_{Z}$ is the Albanese torsor of $Z$, hence $\upic^{1}_{Z}(k) = \emptyset$. This implies that the index of $Z$ is not $1$.
	
	$(3) \Leftrightarrow (4)$ Follows directly from Proposition~\ref{prop:toricbrauer} plus the fact that, by a theorem of Roquette and Lichtenbaum \cite{lichtenbaum}, the relative Brauer group is trivial if and only if the index is one.
\end{proof}

\section{Towers in abelian varieties}\label{sect:tower}

Let $A$ be an abelian variety over a field $k$ finite over $\QQ_{p}$, and $X \subset A$ a closed subvariety. We write
\begin{itemize}
	\item $X - X$ for the schematic image of the morphism $X \times X \to A$ defined by $(x,y) \mapsto y - x$.
	\item $n^{-1} X$ for the inverse image of $X$ along the multiplication by $n$ homomorphism $n: A \to A$.
\end{itemize}

The purpose of this section is to prove the following Theorem~\ref{thm:tower}.

\begin{theorem}\label{thm:tower}
	Let $k$ be a finite extension of $\QQ_{p}$, $A$ an abelian variety over $k$, $X \subset A$ a closed subvariety such that $X \cap A[p'] = \emptyset$. Assume either that
	\begin{itemize}
		\item $A$ has potentially good reduction, or
		\item $X - X$ contains no translates of positive dimensional abelian subvarieties over $\bar{k}$.
	\end{itemize}
	For every $n >> 0$ large enough, $p^{-n}X$ has index divisible by $p$.
\end{theorem}

\begin{remark}
	It is rather easy to show that, since $X \cap A[p'] = \emptyset$, for every closed point $x \in X$ there exists $n$ such that the fiber $p^{-n}\{x\}$ has index divisible by $p$. By compactness, we can find a single $n$ which works uniformly for every $k'$-point, where $k'$ is any fixed finite extension of $k$. The point of Theorem~\ref{thm:tower} is finding a single $n$ which works uniformly for all closed points at the same time.
\end{remark}

Clearly, the index decreases after base change. Hence, we may freely replace $k$ with finite extensions.

Up to extending $k$, we may then assume that $k$ contains a primitive $p$-th root of $1$, and that $A$ has semistable reduction with Néron model $\cA$.

\subsection{Local analysis}\label{sect:local}

First, we study a local analogue of Theorem~\ref{thm:tower}. Let us begin by fixing some notation.

\begin{notation}\label{not:not-loc}
	In subsection \ref{sect:local}, we use the following notation.
	\begin{itemize}
		\item $R$ is a complete DVR of mixed characteristic $(0,p)$ with fraction field $K$, residue field $F$, uniformizing parameter $\pi$, and valuation $v$ normalized by $v(\pi)=1$. We assume that $R$ contains a primitive $p$-th root of $1$.
		\item We write $\CC_{K}$ for the completion of $\bar{K}$ with respect to the unique prolongation of the valuation to $\bar{K}$, and $\fro \subset \CC_{K}$ for the elements of strictly positive valuation. We still write $v$ for the valuation on $\CC_{K}$.
		\item $G \simeq \spf R[[T_{1}, \dots, T_{g}]]$ is a smooth, commutative formal group over $R$.
		\item We shorten the variables $(T_{1}, \dots, T_{g})$ as a variable-vector $T = (T_{1}, \dots, T_{g})$, and similarly $U = (U_{1}, \dots, U_{g})$ for vectors in $\CC_{K}^{\oplus g}$.
		\item Products of vectors are intended component-wise, and not as scalar products. In particular, for every positive integer $m$, we write $T^{m} = (T_{1}^{m}, \dots, T_{g}^{m})$, and similarly for $U$. 
		\item $\Phi(T) =(\Phi_{1}(T), \dots, \Phi_{g}(T))\in R[[T]]^{\oplus g}$ is the formal power series of the multiplication by $p$ map $[p]:G \to G$. 
	\end{itemize}
\end{notation}

We are interested in the case $R=\widehat{\cO_{k^{\rm nr}}}$, the completed maximal unramified valuation ring, with $G$ the formal group at the identity of $\cA_R$.

Loosely speaking, the purpose of this section is to prove that a subvariety of the formal group whose distance from $0$ is larger than $|\pi|^{c}$ cannot contain prime-to-$p$ points which are $p^{2c}$-divisible (see Proposition~\ref{prop:local-bound}), where the distance is intended in the sense of \cite{voloch, scanlon}. 

\begin{remark}\label{rem:counterintuitive}
	Our target is counterintuitive, at least for $G = \hat{\GG}_{m}$. We can find points $U \in G(S) = \frm_{S}^{\oplus g}$ with $|U|$ very close to $1$ and $\deg S/R$ prime to $p$, so that $|[p^{2c}]U|$ is again close to $1$, and consider the subvariety $X$ whose points are the Galois conjugates of $[p^{2c}]U$. For any $x \in X$, we have that $|x|$ is close to $1$. This looks strange, because all the points of $X$ are divisible by $p^{2c}$, which seems to contradict the above.

	This is not the case. The distance of $X$ from $0$ is not $\inf_{x\in X}|x|$, which is indeed large: rather, $d(X,0)=\sup_f|f(0)|$, where $f$ varies among integral functions vanishing on $X$. For instance, in dimension one put $y=[p^{2c}]U$ and let $K(y)$ be its minimal field of definition. The minimal polynomial of $y$ defines $X$, so
	\[d(X,0)=\left|\rN_{K(y)/K}(y)\right|=|y|^{[K(y):K]}.\]
	This can be very small when $[K(y):K]$ is large even if $|y|$ is close to $1$.
\end{remark}

\begin{remark}\label{rem:ai}
	A crucial ingredient of the proof of Theorem~\ref{thm:tower} is a theorem by Scanlon giving a uniform lower bound on the distance from $X$ to torsion points. The correct distance is lower than the ``naive'' distance described in Remark~\ref{rem:counterintuitive}, hence Scanlon's theorem is a priori stronger than its ``naive'' counterpart. We need the strong version to prove Theorem~\ref{thm:tower}, while the weak version does not seem to be sufficient: this hindered progress for a long time. It turns out that the weak version implies the strong version \cite[p. 639]{voloch}, but this requires a theorem of Greenberg of which we were not aware \cite{greenberg}. ChatGPT's main contribution to the paper was to help me realize that I was using the weak version of Scanlon's theorem while I needed the strong one.
\end{remark}

\begin{remark}\label{rem:proofsketch}
	Before starting the proof, let us sketch the idea. Fix a constant $c > 0$. First we prove that, for $n>>0$ large enough, $[p^{n}]: G \to G$ can be approximated modulo $\pi^{c+1}$ by power series in $T^{p^{c}}$. This follows from standard arguments. Hence, the same is true for $f(\Phi^{\circ n}(T))$, where $f \in R[[T]]$ is any power series. It turns out that $n = 2c$ is sufficient.
	
	The idea is that, if $S/R$ has degree prime with $p$ and $s \in \frm_{S}$, we have $\tr(s^{p^{c}}) \in \pi^{c+1}R$. More generally, this is true for every symmetric polynomial with zero constant term evaluated on the Galois conjugates of $s^{p^{c}}$, and this generalized statement has an easy proof by induction on $c$.
	
	Let $f \in R[[T]]$ be a function vanishing on $X$ with $f(0) \neq 0$. By the above, if $U \in G$ is a prime-to-$p$ point and $c = v(f(0))$, we get that 
	\[v(\tr(f([p^{2c}]U))) = v(\tr(f(0))) = v(f(0)).\]
	This implies that $f$ does not vanish on $[p^{2c}]U$, hence $[p^{2c}]U \not \in X$.
\end{remark}

\begin{lemma}\label{lem:frobenius}
	Let $n \ge 0$ be an integer. There exists $\Psi(T) \in R[[T]]^{\oplus g}$ such that
	\[\Phi^{\circ n}(T) \cong \Psi(T^{p^{n}}) \pmod{\pi}.\] 
\end{lemma}

\begin{proof}
	On the special fiber, we may factorize $[p^{n}] = V \circ F$, where $F: G \to G^{(p^{n})}$ and $V: G^{(p^{n})} \to G$ are the Frobenius and Verschiebung homomorphisms \cite[Chapter 2, \S 5]{demazure}. The formal group $G^{(p^{n})}$ can be identified with the formal spectrum of $F[[T^{p^{n}}]]$, with Frobenius simply given by the inclusion $F[[T^{p^{n}}]] \subset F[[T]]$. This implies that the reduction of $\Phi^{\circ n}$ is a series in $T^{p^{n}}$, which is exactly the statement.
\end{proof}

\begin{lemma}\label{lem:pgain}
	Let $a \ge 1$ be an integer. If $U, V \in \fro^{\oplus g}$ are vectors with
	\[U \cong V \pmod{\pi^{a}},\]
	then 
	\[\Phi(U) \cong \Phi(V) \pmod{\pi^{a+1}}.\]
\end{lemma}

\begin{proof}
	By Lemma~\ref{lem:frobenius},
	\[\Phi(T) = \Psi(T^{p}) + \pi \Gamma(T)\]
	for some $\Psi,\Gamma \in R[[T]]^{\oplus g}$.
	
	By hypothesis, $U = V + \pi^{a} W$ for some integral vector $W \in \CC_{K}^{\oplus g}$. Then
	\[\Phi(U) = \Psi((V + \pi^{a} W)^{p}) + \pi \Gamma(V + \pi^{a} W) .\]
	We may write 
	\[(V + \pi^{a} W)^{p} = V^{p} + p\pi^{a}WV^{p-1} + \pi^{2a}W_{0} = V^{p} + \pi^{a+1}W_{1}\]
	for some integral vectors $W_0,W_1\in\cO_{\CC_K}^{\oplus g}$, hence
	\[\Phi(U) = \Psi(V^{p} + \pi^{a+1}W_{1}) + \pi\Gamma(V + \pi^{a}W) \cong \Psi(V^{p}) + \pi\Gamma(V) = \Phi(V) \pmod{\pi^{a+1}}.\]
\end{proof}

\begin{lemma}\label{lem:max-mod}
	Let $\Psi(T) \in R[[T]]$ be a power series, and $a \ge 1$ an integer. Then
	\[\Psi(T) \cong 0 \pmod{\pi^{a}}\]
	if and only if, for every $U \in \fro^{\oplus g}$, we have
	\[\Psi(U) \cong 0 \pmod{\pi^{a}}.\]	
\end{lemma}

\begin{proof}
	Write $\Psi(T) = \sum_{I}a_{I}T^{I}$. For every $0 < r \le 1$, write $D_{r} \subset \CC_{K}^{\oplus g}$ for the closed disk of radius $r$. While $\Psi$ might not converge on $D_{1}$, it does converge on $D_{r}$ for $r < 1$. By \cite[\S 5.1.4, Corollary 6]{bgr}, we then have $\sup_{U \in D_{r}}|\Psi(U)| = \sup_{I} |a_{I}|r^{|I|}$. Passing to the limit for $r \to 1$, we get $\sup_{U \in \fro^{\oplus g}}|\Psi(U)| = \sup_{I} |a_{I}|$. We have $\Psi(T) \cong 0$ if and only if $|a_{I}| \le |\pi^{a}|$ for every $I$, if and only if $|\Psi(U)| \le |\pi^{a}|$ for every $U \in \fro^{\oplus g}$, if and only if $\Psi(U) \cong 0$ for every $U$.
\end{proof}

\begin{lemma}\label{lem:phi-psi}
	Let $n,c \ge 0$ be integers. There exists $\Psi(T) \in R[[T]]^{\oplus g}$ such that
	\[\Phi^{\circ n + c}(T) \cong \Psi(T^{p^{n}}) \pmod{\pi^{c+1}}.\] 
\end{lemma}

\begin{proof}
	By Lemma~\ref{lem:frobenius}, we may find $\Psi_{0}$ such that
	\[\Phi^{\circ n}(T) \cong \Psi_{0}(T^{p^{n}}) \pmod{\pi}.\]
	By Lemma~\ref{lem:max-mod}, this is equivalent to
	\[\Phi^{\circ n}(U) \cong \Psi_{0}(U^{p^{n}}) \pmod{\pi}\]
	for every $U \in \fro^{\oplus g}$. Write $\Psi = \Phi^{\circ c} \circ \Psi_{0}$. By Lemma~\ref{lem:pgain}, we get
	\[\Phi^{\circ n + c}(U) \cong \Phi^{\circ c} \circ \Psi_{0}(U^{p^{n}}) = \Psi(U^{p^{n}}) \pmod{\pi^{c+1}},\]
	for every $U \in \fro^{\oplus g}$, hence
	\[\Phi^{\circ n + c}(T) \cong \Psi(T^{p^{n}})\pmod{\pi^{c+1}}.\]
\end{proof}

\begin{lemma}\label{lem:symmetric}
	Let $S/R$ be an extension of degree prime with $p$, and $s \in \frm_{S}$ an element in the maximal ideal of $S$. Fix an embedding $S \subset \CC_{K}$, and let $s_{1}, \dots, s_{e} \in \CC_{K}$ be the Galois conjugates of $s$, with $s_{i} \neq s_{j}$ if $i \neq j$. 
	
	Let $a \ge 1$ be an integer. Assume that, for every symmetric polynomial $q \in \ZZ[t_{1}, \dots, t_{e}]$ with $q(0, \dots, 0) = 0$, we have $q(s_{1}, \dots, s_{e}) \in \pi^{a} R$. 
	
	Then $s_{1}^{p}, \dots, s_{e}^{p}$ are the Galois conjugates of $s^{p}$, with $s_{i}^{p} \neq s_{j}^{p}$ if $i \neq j$, and for every symmetric polynomial $q$ with $q(0, \dots, 0) = 0$ we have $q(s_{1}^{p}, \dots, s_{e}^{p}) \in \pi^{a+1}R$.
\end{lemma}

\begin{proof}
	Since $R$ contains a primitive $p$-th root of $1$, then $K(s_{i})/K(s_{i}^{p})$ is Galois of degree either $1$ or $p$: since the degree is prime with $p$, we must have $K(s_{i}^{p}) = K(s_{i})$ for every $i$. This implies that $s_{i}^{p} \neq s_{j}^{p}$ if $i \neq j$, because otherwise a Galois automorphism sending $s_{i}$ to $s_{j}$ fixes $K(s_{i}^{p}) = K(s_{i})$, which is absurd. Since $s^{p}$ has $[K(s^{p}):K] = [K(s):K] = e$ conjugates, we get that $s_{1}^{p}, \dots, s_{e}^{p}$ are the conjugates of $s^{p}$, with $s_{i}^{p} \neq s_{j}^{p}$ if $i \neq j$.
	
	It is enough to prove the statement for a symmetric polynomial $q$ all of whose coefficients are $1$, because every symmetric polynomial can then be expressed as a linear combination of these ones. For such a $q$, we have
	\[q(s_{1}, \dots, s_{e})^{p} = q(s_{1}^{p}, \dots, s_{e}^{p}) + ph(s_{1},\dots,s_{e}),\]
	where $h$ is some symmetric polynomial with $h(0, \dots, 0) = 0$. This implies that $q(s_{1}^{p}, \dots, s_{e}^{p}) \in \pi^{a+1}R$, hence we conclude.
\end{proof}

\begin{corollary}\label{cor:trace}
	Let $a\ge0$ be an integer. Let $S/R$ be an extension of degree prime with $p$, and $s \in \frm_{S}$ an element in the maximal ideal of $S$. Then
	\[\tr_{S/R}\left(s^{p^{a}}\right) \in \pi^{a+1}R.\]	
\end{corollary}

\begin{proof}
	Let $s_{1}, \dots, s_{e} \in \CC_{K}$ be the Galois conjugates of $s$ as in Lemma~\ref{lem:symmetric}. Notice that $s_{i}$ has positive valuation for every $i$, because $R$ is complete and hence the unique extension of the valuation to $\bar{K}$ is Galois invariant.
	
	This implies that $s$ satisfies the hypothesis of Lemma~\ref{lem:symmetric} with $a = 1$. By induction on $a$, using Lemma~\ref{lem:symmetric} we get that $f(s_{1}^{p^{a}}, \dots, s_{e}^{p^{a}}) \in \pi^{a+1}R$ for every symmetric polynomial $f \in \ZZ[t_{1}, \dots, t_{e}]$ with $f(0, \dots, 0) = 0$. In particular, $\tr_{S/R}(s^{p^{a}}) \in \pi^{a+1}R$. Notice that $\tr_{S/R}(s^{p^{a}})$ is not necessarily the sum of the Galois conjugates of $s^{p^{a}}$ but a multiple of it, e.g. if $s^{p^{a}} \in R$. Still, it is a symmetric polynomial of the conjugates with zero constant term, and that's enough. 
\end{proof}

\begin{proposition}\label{prop:local-bound}
	Let $c \ge 0$ be an integer, and $f \in R[[T]]$ a power series with $v(f(0, \dots, 0)) \le c$. For every extension $S/R$ of degree $e$ prime with $p$ and every $U \in \frm_{S}^{\oplus g}$, we have 
	\[\tr_{S/R}(f(\Phi^{\circ 2c}(U))) \cong ef(0, \dots, 0) \not \cong 0 \pmod{\pi^{c+1}}.\]
	In particular, if $b \in G(S) = \frm_{S}^{\oplus g}$ is a point defined over $S$, then
	\[f([p^{2c}]b) \neq 0.\]
\end{proposition}

\begin{proof}
	Thanks to Lemma~\ref{lem:phi-psi}, we can write
	\[f(\Phi^{\circ 2c}(T)) \cong h(T^{p^{c}}) \pmod{\pi^{c+1}}\]
	for some power series $h(T) \in R[[T]]$. Let us study $h(U^{p^{c}})$.

	The difference $h(U^{p^{c}}) - h(0, \dots, 0)$ is a convergent, infinite sum of terms of the form $\alpha s^{p^{c}}$ for some $\alpha \in R$ and $s \in \frm_{S}$. By Corollary~\ref{cor:trace},
	\[\tr_{S/R}\left(\alpha s^{p^{c}}\right) = \alpha \tr_{S/R}\left(s^{p^{c}}\right) \in \pi^{c+1}R.\]
	
	This implies
	\[\tr_{S/R}(f(\Phi^{\circ 2c}(U))) \cong \tr_{S/R}\left(h(U^{p^{c}})\right) \cong \tr_{S/R}\left(h(0, \dots, 0)\right) = \]
	\[ = e h(0, \dots, 0) \cong e f(\Phi^{\circ 2c}(0, \dots, 0)) = ef(0, \dots, 0) \not \cong 0 \pmod{\pi^{c+1}}.\]
\end{proof}

\subsection{Special fiber}

Write $F$ for the residue field of $k$, $F^{(p')}$ for its maximal prime-to-$p$ extension, and $\cA$ for the Néron model.

\begin{proposition}\label{prop:finiteshallow}
	The group $\cA[p^{\infty}](F^{(p')})$ is finite. As a consequence, there exists an integer $r > 0$ such that $p^{r}\cA(F^{(p')})$ is $p'$-torsion.
\end{proposition}

\begin{proof}
	Let $B \subset \cA_{F}$ be the connected component of the identity, it is enough to prove that $B[p^{\infty}](F^{(p')})$ is finite since the number of connected components is finite.
	
	Let $r$ be the rank of $B[p](\bar{F})$. For every $n$, the group scheme $B[p^{n}]$ gives an action $\gal(\bar{F}/F) \to \GL(r,\ZZ/p^{n})$. Notice that the kernel of $\GL(r,\ZZ/p^{n}) \to \GL(r,\ZZ/p)$ is a $p$-group: a matrix in the kernel has the form $\operatorname{Id} + pM$, and $(\operatorname{Id} + pM)^{p^{m}} = \operatorname{Id}$ for $m$ large enough. Since $\GL(r,\ZZ/p)$ is finite, there exists a finite extension $E/F$ such that the composition $\gal(\bar{E}/E) \to \GL(r,\ZZ/p)$ is trivial. In particular, the image of $\gal(\bar{E}/E) \to \GL(r,\ZZ/p^{n})$ is a $p$-group for every $n$.
		
	This implies that the images of $\gal(\bar{E}/E) \to \GL(r,\ZZ/p^{n})$ and $\gal(\bar{E}/E^{(p)'}) \to \GL(r,\ZZ/p^{n})$ coincide for every $n$. We get that
	\[B[p^{\infty}](F^{(p')}) \subseteq B[p^{\infty}](E^{(p')}) = B[p^{\infty}](E)\]
	is finite.
	
	Hence, we have $p^{r}\cA[p^{\infty}](F^{(p')}) = \{0\}$ for $r$ large enough. Since $\cA(F^{(p')}) = \cA[p^{\infty}](F^{(p')}) \times \cA[p'](F^{(p')})$, we conclude that $p^{r}\cA(F^{(p')}) \subseteq \cA[p'](F^{(p')})$ is $p'$-torsion.
\end{proof}

For every extension $L/k$, denote by $A^{1}(L) \subset \cA(\cO_{L}) \subset A(L)$ the subgroup of integral points which specialize to the identity. 

\begin{proposition}\label{prop:special}
	There exists a positive integer $r$ such that, for any finite extension $L/k$ of degree prime with $p$, the group
	\[p^{r}(A(L)/A^{1}(L))\]
	is $p'$-torsion.
\end{proposition}

\begin{proof}
	Let us first prove that $p^{r}(\cA(\cO_{L})/A^{1}(L))$ is $p'$-torsion for a fixed $r$ not depending on $L$.
	
	We have $\cA(\cO_{L})/A^{1}(L) = \cA_{F}(F')$, where $F'$ is the residue field of $L$. Let $F^{(p')}$ be the maximal prime-to-$p$ algebraic extension of $F$, i.e. the fixed field of the unique $p$-Sylow of $\gal(\bar{F}/F)$. We have $\cA_{F}(F') \subset \cA_{F}(F^{(p')})$, hence we may conclude by applying Proposition~\ref{prop:finiteshallow}.
	
	Now recall that $A$ has semi-stable reduction. Let $\cB$ be the Néron model of $A_{L}$, so that $A(L^{\rm nr}) = \cB(\cO_{L^{\rm nr}})$. By \cite[Theorem 5.3.3.9]{halle-nicaise}, $\cB(\cO_{L^{\rm nr}})/\cA(\cO_{L^{\rm nr}})$ is $e$-torsion, where $e = [L^{\rm nr}:k^{\rm nr}]$ is prime with $p$. Hence, $eA(L) \subset \cA(\cO_{L})$. By the above, $ep^{r}(A(L)/A^{1}(L)) \subseteq p^{r}(\cA(\cO_{L})/A^{1}(L))$ is $p'$-torsion, which implies that $p^{r}(A(L)/A^{1}(L))$ is $p'$-torsion as well.
\end{proof}

\begin{corollary}\label{cor:close-to-torsion}
	Notation as in Proposition~\ref{prop:special}. For every $a \in A(L)$, there exists $\tau \in A[p'](L)$ and $b \in A^{1}(L)$ such that
	\[p^{r}a = \tau + b.\]
\end{corollary}

\begin{proof}
	By Proposition~\ref{prop:special}, there exists an integer $m$ prime with $p$ such that $m p^{r} a = b_{0} \in A^{1}(L)$. Since $m$ is prime with $p$, multiplication by $m$ is invertible on the formal group $A^{1}(L) \simeq (\spf \cO_{L}[[T_{1}, \dots, T_{g}]])(\cO_{L})$ (this follows from the formal inverse function theorem, see for instance \cite[Appendix A.4]{hazewinkel}), hence there exists a unique $b \in A^{1}(L)$ such that $mb = b_{0}$. It follows that $\tau = p^{r}a - b$ is $m$-torsion, because $m\tau = mp^{r}a - mb = b_{0} - b_{0} = 0$.
\end{proof}

\subsection{Torsion capture}

Write $R$ for the completion of $\cO_{k^{\rm nr}}$, where $k^{\rm nr}$ is the maximally unramified extension, and $K=\operatorname{Frac}(R)$.

Let us recall a theorem of T. Scanlon \cite{scanlon}.

\begin{theorem}[{T. Scanlon}]\label{thm:scanlon}
	Let $Y \subset A$ be a subvariety. There exists a uniform lower bound on the distance between $Y$ and points of $\cA[p'](R) \setminus Y$. More explicitly, there exists an integer $c(Y) > 0$ such that, for every $\tau \in \cA[p'](R)$, either $\tau \in Y$ or there exists an integral function $f$ on a neighbourhood of $\tau$ in $\cA$ which vanishes on $Y$ and satisfies $v(f(\tau)) \le c(Y)$.
\end{theorem}

We call $c(Y)$ \emph{Scanlon's constant} for $Y$. Scanlon's theorem gives us a way to ``capture'' $p'$-torsion points defined over $R$.

\begin{corollary}\label{cor:torsion-capture}
	Notation as in Theorem~\ref{thm:scanlon}. Let $\tau \in \cA[p'](R)$ an integral $p'$-torsion point over $R$, $L/K$ be a finite extension of degree prime to $p$, $b \in A^{1}(L)$ an element in the local formal group. If $\tau + p^{2c(Y)}b \in Y$, then $\tau \in Y$. 
\end{corollary}

\begin{proof}
	By Theorem~\ref{thm:scanlon}, if $\tau \not \in Y$ there exists an integral function $f_{0}$ on a neighbourhood $U_{0}$ of $\tau$ in $\cA$ which vanishes on $Y$ and satisfies $v(f_{0}(\tau)) \le c(Y)$. Consider the neighbourhood $U = U_{0} - \tau$ of $0$, and the integral function $f(u) = f_{0}(\tau + u)$ on $U$; in particular, $v(f(0)) \le c(Y)$. By Proposition~\ref{prop:local-bound}, $f(p^{2c(Y)}b) \neq 0$, which is absurd because $\tau + p^{2c(Y)}b \in Y$ and hence $f(p^{2c(Y)}b) = f_{0}(\tau + p^{2c(Y)}b) = 0$.
\end{proof}

\subsection{Good reduction}\label{sect:integral}

Assume that $A$ has good reduction, so that every point is integral and prime-to-$p$ torsion is defined over $R=\widehat{\cO_{k^{\rm nr}}}$. Let $r$ be the integer given by Proposition~\ref{prop:special}, and $c = c(X)$. 

Choose $L/k$ any finite extension of degree prime with $p$, $a \in A(L)$ a point. By Corollary~\ref{cor:close-to-torsion}, $p^{r}a = \tau + b$ with $\tau \in A[p'](L)$ and $b \in A^{1}(L)$. If $p^{r+2c}a = p^{2c}\tau + p^{2c}b \in X$, then $p^{2c}\tau \in X$ by Corollary~\ref{cor:torsion-capture}, which contradicts the hypothesis that $X \cap A[p'] = \emptyset$. This concludes the proof of the first half of Theorem~\ref{thm:tower}.

\subsection{$(X-X)_{\bar{k}}$ contains not translates of abelian subvarieties}

Let us prove Theorem~\ref{thm:tower} under the assumption that $(X-X)_{\bar{k}}$ contains no translate of a positive dimensional abelian subvariety. By Manin--Mumford \cite[p.~327]{raynaud}, $X-X$ contains finitely many torsion points. Fix $m$ with $p \nmid m$ such that, if $\tau \in (X - X) \cap A[p']$, then $m\tau = 0$. 

Let $r$ be the integer given by Proposition~\ref{prop:special}, and $c = \max\{c(mX),c(X-X)\}$ a Scanlon's constant working both for $mX \subset A$ and $X - X$. 

Let us prove that, if $L/k$ is a finite extension with $p \nmid [L:k]$ and $a \in A(L)$ is a point, then $p^{r+2c}a \not \in X$; equivalently, $p^{-r-2c}X$ has index divisible by $p$. 

By Corollary~\ref{cor:close-to-torsion}, there exists $\tau_{0} \in A[p'](L)$ and $b \in A^{1}(L)$ such that $p^{r}a = \tau_{0} + b$. Write $\tau = p^{2c}\tau_{0}$, we have $p^{r+2c}a = \tau + p^{2c}b \in X$.

\begin{lemma}\label{lem:inertia-torsion}
	If $\tau \in A[p'](\bar{k})$ and $\alpha, \beta \in \gal(\bar{k}/k^{\rm nr})$, then $\alpha (\beta \tau - \tau) = \beta \tau - \tau$.
\end{lemma}

\begin{proof}
	Clearly, we may reduce to the case in which $\tau \in A[\ell^{d}]$ for some prime $\ell \neq p$, which in turn reduces to the analogous statement on $T_{\ell}A$. By Grothendieck's characterization of abelian varieties with semistable reduction \cite[Exposé IX, Proposition 3.5]{sga7.1}, we have a short exact sequence
	\[0 \to N \to T_{\ell}A \to M \to 0\]
	of Galois modules where the Galois actions of the inertia on $N$ and $M$ are both trivial. If $v \in T_{\ell}A$, we get that $\beta v - v \in N$, and hence $\alpha(\beta v - v) - (\beta - v) = 0$.
\end{proof}

Let us work in $L^{\rm nr}$. Notice that $L^{\rm nr}/k^{\rm nr}$ is a cyclic Galois extension of order prime with $p$. Choose $\gamma \in \gal(L^{\rm nr}/k^{\rm nr})$ a generator.

By Lemma~\ref{lem:inertia-torsion}, $\delta = \gamma \tau - \tau$ is Galois invariant, hence $\delta \in A[p'](k^{\rm nr}) = \cA[p'](R)$. Notice that 
\[\gamma p^{r+2c}a - p^{r+2c}a = \delta + p^{2c}(\gamma b - b) \in X - X.\]
By Corollary~\ref{cor:torsion-capture}, $\delta \in X - X$. This implies that 
\[\gamma m\tau - m \tau = m\delta = 0,\]
hence $m\tau$ is Galois invariant because $\gamma$ is a generator, i.e. $m\tau \in A[p'](k^{\rm nr}) = \cA[p'](R)$. 

Clearly, we have $mp^{r+2c}a = m\tau + p^{2c}mb \in mX$, with $m\tau$ integral. By Corollary~\ref{cor:torsion-capture}, we get that $m\tau \in mX$, which is a contradiction because $X \cap A[p'] = \emptyset$ implies $mX \cap A[p'] = \emptyset$.


\section{Curves over $p$-adic fields}\label{sect:curves}

Let $X$ be a hyperbolic curve over a field $k$ finite over $\QQ_{p}$. The goal of this section is to prove Theorem~\ref{thm:torQp}, i.e. that $X(k) = \cS^{\tor}_{X/k}$. At the same time, we are going to prove the first assertion of Theorem~\ref{thm:nbd-p}, i.e. that every non-geometric Galois section $s \in \cS^{\et}_{X/k}$ has an étale neighbourhood $Y \to X$ such that $p$ divides the index of $Y$. By Theorem~\ref{thm:toretinj}, Theorem~\ref{thm:stix} and Proposition~\ref{prop:torindex}, it is immediate to see that Theorem~\ref{thm:torQp} is equivalent to the first assertion of Theorem~\ref{thm:nbd-p}. 

\subsection{$p$-neighbourhoods}

Furthermore, we will prove the second half of Theorem~\ref{thm:nbd-p}. Namely, if $X(k) = \emptyset$ we will find a $p$-neighbourhood $Y \to X$ of index divisible by $p$, in the following sense. 

\begin{definition}
	A $p$-neighbourhood $(Y,r)$ is a neighbourhood such that $\deg(Y/X)$ is a power of $p$ and $Y \to X$ is geometrically Galois. 
\end{definition}

Notice that, if $Y \to X$ is a $p$-neighbourhood and $Y' \to Y$ is a further $p$-neighbourhood, there exists a geometric Galois closure whose degree is a power of $p$. This implies the existence of a further $p$-neighbourhood $Y'' \to Y'$ such that $Y'' \to X$ is a $p$-neighbourhood. This means that we can freely replace $(X,s)$ with a $p$-neighbourhood.

Let $\pi_{1}^{(p)}(X) \to \gal(\bar{k}/k)$ be the pushout of the maximal pro-$p$ quotient $\pi_{1}^{(p)}(X_{\bar{k}})$ of $\pi_{1}(X_{\bar{k}})$, then we can think of $p$-neighbourhoods as the neighbourhoods of the induced section $s^{(p)}: \gal(\bar{k}/k) \to \pi_{1}^{(p)}(X)$.

\subsection{Abelian neighbourhoods}

If $s \in \cS^{\tor}_{X/k}$ is a toric Galois section, for every étale neighbourhood $(Y,r)$ the image of $r$ in $\cS^{\tor}_{\upic^{1}_{Y}/k}$ is associated with a unique rational point in $\upic^{1}_{Y}$ by Theorem~\ref{thm:torab+}.

\begin{convention}
    Given a toric Galois section $s$ of a curve $X$, for every étale neighbourhood $(Y,r)$ we write $J_{Y} = \upic^{1}_{Y}$ and consider it as an abelian variety where the origin is the unique rational point associated with the image of $r$ in $\cS^{\et}_{J_{Y}/k}$.
\end{convention}

We are going to use the following observation repeatedly.

\begin{lemma}\label{lem:nbd}
	Let $s \in \cS^{\tor}_{X/k}$ be a toric Galois section of a curve $X$ and $n \ge 0$ an integer, consider the multiplication by $p^{n}$ map $p^{n} : J_{X} \to J_{X}$ and the inverse image $p^{-n}X \subset J_{X}$. The map $p^{-n}X \to X$ is a $p$-neighbourhood of $s$.
\end{lemma}

\begin{proof}
	This follows from the fact that $\pi_{1}(p^{-n}X) \subset \pi_{1}(X)$ is the inverse image of $\pi_{1}(J_{X}) \subset \pi_{1}(J_{X})$ (where the latter embedding is given by $p^{n}: J_{X} \to J_{X}$). In fact, the image of $s$ in $\cS^{\et}_{J_{X}/k}$ is by convention geometric associated with the origin, and the origin clearly lifts along $p^{n}: J_{X} \to J_{X}$.
\end{proof}

\subsection{Avoiding $p'$-torsion}

Given an abelian variety $A$, we denote by $A[p']$ the subgroup of torsion points of order prime with $p$.

\begin{proposition}\label{prop:avoidtorsion}
    Let $X$ be a hyperbolic curve of a field $k$ finite over $\QQ_{p}$, and $s \in \cS_{X/k}^{\tor}$ a toric Galois section. If $s$ is not geometric, there exists an étale neighbourhood $(Y,r)$ of $s$ such that 
    \[Y\cap J_{Y}[p']=\emptyset.\]
    Furthermore, any further étale neighbourhood of $(Y,r)$ avoids $p'$-torsion.
    
    If $X(k) = \emptyset$, we can choose $Y \to X$ a $p$-neighbourhood. 
\end{proposition}

\begin{proof}
	Notice that if $f:Z \to Y$ is a further neighbourhood, then $f(Z \cap J_{Z}[p']) \subset Y \cap J_{Y}[p'] = \emptyset$. Let us prove the existence of $(Y,r)$.

    Since $s$ is not geometric, up to passing to an étale neighbourhood we can reduce to the case in which $X(k) = \emptyset$ \cite[Proposition 2.8 (iv)]{tamagawa} \cite[Corollary 101]{stix}.
    
    By Raynaud's theorem (Manin--Mumford conjecture) \cite[p.~327]{raynaud}, $X\cap J_{X}[p']$ is finite. Hence, to prove the statement is enough to prove that, for every $x\in X\cap J_{X}[p'](\bar{k})$, there exists a $p$-neighbourhood $(Y, r)$ such that $x$ is not in the image of
    \[Y\cap J_{Y}[p'](\bar{k}) \to X\cap J_{X}[p'](\bar{k}).\]
	
    Assume by contradiction that this is false, i.e. $x$ is the image of a point $y\in Y\cap J_{Y}[p'](\bar{k})$ for every $p$-neighbourhood $(Y,r)$. Since $X(k)=\emptyset$, there exists a geometric point $x'\neq x$ which is a Galois conjugate of $x$; in particular, $x'$ satisfies the same property as $x$, i.e. it is the image of a point $y'\in Y\cap J_{Y}[p'](\bar{k})$ for every $p$-neighbourhood $(Y,r)$. Notice that, by construction, $[y']-[y]\in\upic^{0}_{Y}[p']$.
    
    For every Galois $p$-cover $Z\to X_{\bar{k}}$ of $X_{\bar{k}}$, there exists a $p$-neighbourhood $(Y,r)$ of $s$ with a factorization $Y_{\bar{k}}\to Z\to X_{\bar{k}}$; this follows from the fact that $s^{(p)}(\gal(\bar{k}/k)) \subset \pi_{1}^{(p)}(X)$ is the intersection of the fundamental groups of the $p$-neighbourhoods of $s$. This is proved in \cite[\S 4.4]{stix} for $\pi_{1}(X)$; the proof of $\pi_{1}^{(p)}(X)$ is identical. 
    
    Thanks to what we have said above, the two geometric points $x,x'$ enjoy a particular property: for every Galois $p$-cover $Z\to X_{\bar{k}}$, there are two points $z,z'\in Z$ mapping to $x,x'$ and such that $[z']-[z]\in\upic^{0}_{Z}[p']$. This is in contradiction with the following Proposition~\ref{prop:ettors}.
\end{proof}

\begin{proposition}\label{prop:ettors}
	Let $X$ be a hyperbolic curve over a field $k$ of characteristic $0$, and $p$ a prime number. Let $x,x'\in X(\bar{k})$ be two distinct geometric points. 
	
	There exists a Galois $p$-cover $Z\to X_{\bar{k}}$ such that, for every pair of geometric points $z,z'\in Z$ over $x,x'$, the difference $[z']-[z]\in\upic^{0}_{Z}$ is not in $\upic^{0}_{Z}[p']$. Furthermore, we may assume that $Z\to X_{\bar{k}}$ is a Galois covering whose degree is a power of $p$.
\end{proposition}

\begin{proof}
	We can reduce to the case in which $k$ is a finite extension of $\QQ_{p}$: in fact, everything is defined over a field of finite type over $\QQ$, and such a field can be embedded in a finite extension of $\QQ_{p}$. Furthermore, up to a finite extension of $k$ we might assume that $x,x'$ are $k$-rational.
	
	We use $x$ as the base point. We have a short exact sequence
	\[1\to \pi_{1}^{p}(X_{\bar{k}})\to \pi_{1}^{p}(X)\to\gal(\bar{k}/k)\to 1,\]
	a space of pro-$p$ Galois sections $\cS^{\et,p}_{X/k}$ and a section map $X(k)\to\cS^{\et,p}_{X/k}$. S. Mochizuki has proved that this pro-$p$ section map is injective under our hypothesis \cite[Theorem 19.1]{mochizuki}.
	
	Write $(X_{0},x_{0})=(X,x)$ and define by recursion a pointed finite étale cover $(X_{i+1},x_{i+1})\to (X_{i},x_{i})$ as follows: we consider the embedding $X_{i}\subset \upic^{1}_{X_{i}}$, and define $X_{i+1}\to X_{i}$ as the pullback of the multiplication by $p$ map $p:\upic^{1}_{X_{i}}\to \upic^{1}_{X_{i}}$ with respect to the base point $x_{i}$. The point $x_{i+1}$ over $x_{i}$ is the one associated with the origin of $\upic^{1}_{X_{i}}$. Write $(X_{\infty},x_{\infty})$ for the projective limit $\projlim_{i}(X_{i},x_{i})$. Notice that $X_{i} \to X$ is a geometrically Galois $p$-cover, and that every Galois $p$-cover of $X_{\bar{k}}$ is dominated by some $X_{i,\bar{k}}$: this follows from the fact that $p$-groups are nilpotent.
	
	It follows that the base change $X_{\infty,\bar{k}}$ is a connected pro-finite étale Galois $p$-cover of $X_{\bar{k}}$ with Galois group $\pi_{1}^{p}(X_{\bar{k}},x)$. So to speak, $X_{\infty}$ the pro-$p$ decomposition tower cf. \cite[\S 4.4]{stix} of the pro-$p$ Galois section associated with $x\in X(k)$.
	
	Assume by contradiction that the statement is false, i.e. for every finite étale Galois $p$-cover $Z\to X_{\bar{k}}$ there exists a pair of geometric points $z,z'\in Z$ over $x,x'$ such that $[z]-[z']\in\upic^{0}_{Z}[p']$. Consider the case $Z=X_{1,\bar{k}}$, since the cover is Galois we may assume that $z=x_{1}$, write $x_{1}'=z'$. We have a cartesian diagram
	\[\begin{tikzcd}
		X_{1}\rar["j_{1}"]\dar	&	\operatorname{Pic}^{1}_{X_{0}}\dar["p"]	\\
		X_{0}\rar				&	\operatorname{Pic}^{1}_{X_{0}}
	\end{tikzcd}\]
	where we denote by $j_{1}:X_{1}\to \upic^{1}_{X_{0}}$ the upper horizontal arrow. If $a\neq x_{1}'\in X_{1}$ is a different point over $x_{0}'$, then $j_{1}(a)-j_{1}(x_{1}')$ is a \emph{non-trivial} $p$-torsion element of $\upic^{1}_{X_{0}}$, and hence $[a]-[x_{1}']\notin \upic^{0}_{X_{1}}[p']$. This implies that 
	\[[a]-[x_{1}]=([a]-[x_{1}'])+([x_{1}']-[x_{1}]) \notin \upic^{0}_{X_{1}}[p'],\]
	since $[x_{1}']-[x_{1}]\in\upic^{0}_{X_{1}}[p']$. It follows that $x_{1}'$ is the \emph{unique} point over $x_{0}'$ such that $[x_{1}']-[x_{1}]\in\upic^{0}_{X_{1}}[p']$. Since $x_{1}$ is rational and hence Galois invariant, we get that $x_{1}'$ is Galois invariant and hence rational as well. 
	
	By repeating the process, we obtain a tower $(x_{i}')_{i}$ of rational points of $(X_{i})_{i}$ inducing a rational point $x_{\infty}'\in X_{\infty}(k)$ over $x'$. Analogously to what happens for classical Galois sections, this implies that $x$ and $x'$ induce the same pro-$p$ Galois section in $\cS^{\et,p}_{X/k}$, see e.g. \cite[Lemma 52]{stix}. In fact, the pro-$p$ Galois sections of $x$ and $x'$ both factorize through the image $H\subset \pi_{1}^{p}(X)$ of $\pi_{1}(X_{\infty})$ inside $\pi_{1}^{p}(X)$, and $H$ maps bijectively on $\gal(\bar{k}/k)$ since $X_{\infty,\bar{k}}$ is the universal pro-$p$ cover of $X_{\bar{k}}$. This is in contradiction with Mochizuki's theorem \cite[Theorem 19.1]{mochizuki}.
\end{proof}

Theorem~\ref{thm:etale-sep} follows directly from Proposition~\ref{prop:ettors}.

\subsection{Theorem~\ref{thm:torQp} and the first half of Theorem~\ref{thm:nbd-p}}

Thanks to Proposition~\ref{prop:avoidtorsion}, we may assume that $X \cap J_{X}[p'] = \emptyset$. 

\begin{lemma}\label{lem:not-hyp}
	Let $C$ be a smooth projective curve over an algebraically closed field of characteristic $0$. If $H \to C$ is an étale Galois cover of degree $d$ with $H$ hyperelliptic of genus at least $2$, then $d \mid 4$.
\end{lemma}

\begin{proof}
	Write $g_{H}$, $g_{C}$ for the genera. Since the cover is étale, its Galois group acts freely on the $2g_{H} + 2$ points fixed by the hyperelliptic involution, hence $d \mid 2g_{H} + 2$. By Riemann--Hurwitz,
	\[d(2g_{C} - 2) = 2g_{H} - 2 = 2g_{H} + 2 - 4,\]
	hence $d \mid 4$.
\end{proof}

The cover $p^{-2}X \to X$ is étale, geometrically Galois, its degree $p^{4g}$ does not divide $4$, and it is an étale neighbourhood of $s$. Furthermore, by Proposition~\ref{prop:avoidtorsion} we have $p^{-2}X \cap J_{p^{-2}X}[p'] = \emptyset$ as well. Hence, up to replacing $X$ with $p^{-2}X$, we might assume that $X$ is not hyperelliptic. 

Since $X$ is not hyperelliptic, the natural map $X \times X \to X - X \subset J_{X}$ identifies $X - X \setminus \{0\}$ with $X \times X \setminus \Delta$, where $\Delta$ is the diagonal. In particular, $(X-X)_{\bar{k}}$ cannot contain genus $1$ curves, because such a curve would dominate $X$ which has genus $> 1$. Furthermore, $X-X$ is not an abelian surface, because it generates $J_{X}$ and $\dim J_{X} = g > 2$ (again because $X$ is not hyperelliptic). 

Since $\dim X-X = 2$, we get that $(X-X)_{\bar{k}}$ does not contain translates of positive dimensional abelian varieties. By Theorem~\ref{thm:tower}, there exists $n$ such that $p$ divides the index of $p^{-n}X$. Since $p^{-n}X \to X$ is a $p$-neighbourhood, we get a contradiction with Proposition~\ref{prop:torindex}. This concludes the proof of Theorem~\ref{thm:torQp}, and the first half of Theorem~\ref{thm:nbd-p}.

\subsection{Second half of Theorem~\ref{thm:nbd-p}}

We now want to show how to modify the above in order to prove the second half of Theorem~\ref{thm:nbd-p}.

Consider the image $s^{(p)} \in \cS^{\et,(p)}_{X/k}$ of a section $s \in \cS^{\et}_{X/k}$. In the above, we worked under the assumption that $s$ was not geometric. Let us now assume that $s^{(p)}$ is not geometric instead; this is obviously guaranteed if $X(k) = \emptyset$.

In the above, we have also worked under the assumption that $s \in \cS^{\tor}_{X/k} \subseteq \cS^{\et}_{X/k}$. By Proposition~\ref{prop:torindex}, this is equivalent to the fact that, for every neighbourhood $(Y,r)$, the abelianization $r^{\rm ab} \in \cS^{\et}_{\upic^{1}_{Y/k}/k}$ is geometric. Let us work under the weaker assumption that this is true only for $p$-neighbourhoods $(Y,r^{(p)})$ of $s^{(p)}$, so that the abelianization $r^{(p),\rm ab} \in \cS^{\et,(p)}_{\upic^{1}_{Y/k}/k}$ is assumed geometric.

We want to show that $s^{(p)}$ non-geometric is in contradiction with $r^{(p),\rm ab}$ geometric for every $p$-neighbourhood $(Y,r^{(p)})$ of $s^{(p)}$; this proves the second half of Theorem~\ref{thm:nbd-p}.

If $A$ is an abelian variety over $k$, the map $A(k) \to \cS^{\rm et,(p)}_{A/k}$ is not guaranteed to be injective, because all the points of $A[p'](k)$ map to the same Galois $p$-section. Hence, if $(Y,r^{(p)})$ is a $p$-neighbourhood of $s^{(p)}$, there is no natural choice of origin for $\upic^{1}_{Y/k}$, even if $r^{(p), \rm ab}$ is geometric.

However, this is not a problem. Choose any point $a$ of $\upic^{1}_{Y/k}(k)$ associated with $r^{(p), \rm ab}$ as origin, so that we have an abelian variety $J_{Y}$. For every $n$, we have the multiplication by $p^{n}$ homomorphism $J_{Y} \to J_{Y}$. 

Let $b \in \upic^{1}_{Y/k}$ be another choice, and $J_{Y}' \to J_{Y}'$ the corresponding $p^{n}$-homomorphism. We have that $a-b \in \upic^{0}_{Y/k}[p'](k)$ is divisible by $p^{n}$, because $a$, $b$ map to the same $p$-section in $\cS^{\et, (p)}_{\upic^{1}_{Y/k}/k}$. This implies that we can find a commutative diagram
\[\begin{tikzcd}
	J_{Y}\rar["\sim"]\dar["p^{n}"]	&	J_{Y}'\dar["p^{n}"]	\\
	J_{Y}\rar["\sim"]				&	J_{Y}'
\end{tikzcd}\]
where the horizontal arrows are isomorphisms induced by translation on $\upic^{1}_{Y/k}$. In the proof of Theorem~\ref{thm:nbd-p}, we only use the origin of $\upic_{Y/k}^{1}$ to define multiplication by $p^{n}$. Since multiplication by $p^{n}$ is well-defined as long as $r^{(p), \rm ab}$ is geometric for every $p$-neighbourhood, we can repeat the argument and obtain the desired contradiction.

\section{Results over number fields}\label{sect:nf}

Let us now prove Theorem~\ref{thm:selmer}, i.e. that over a number field $k$ we have $\cS_{X/k}^{\tor} = \operatorname{Sel}_{X/k} \subset \cS^{\et}_{X/k}$. The passage from local fields to number fields is delicate since the toric fundamental gerbe does not behave well for arbitrary field extensions which are not algebraic and separable, see Example~\ref{ex:badchange}

Following Koenigsmann \cite{koenigsmann}, we overcome this problem using $p$-adically closed fields.

\begin{lemma}\label{lem:elemdes}
	Let $k$ be a finite extension of $\QQ_{p}$, and $\bar{\QQ}^{k}$ the algebraic closure of $\QQ$ in $k$. Let $X,Y$ be projective varieties over $\bar{\QQ}^{k}$. If there exists a morphism $X_{k}\to Y_{k}$, then there exists a morphism $X\to Y$.
\end{lemma}

\begin{proof}
	By choosing equations for the given morphism $X_k\to Y_k$, the existence of a morphism $X\to Y$ with equations of the same form can be phrased as a sentence in the elementary theory of the field $\bar{\QQ}^{k}$. The inclusion $\bar{\QQ}^{k}\hookrightarrow k$ is elementary \cite[Theorems E.2, E.3]{pop-galois}, hence $X\to Y$ exists if and only if $X_k\to Y_k$ exists.
\end{proof}

\begin{proof}[{Proof of Theorem~\ref{thm:selmer}}]
	Consider a Selmer section $s \in \operatorname{Sel}_{X/k}$, we want to show that $s \in \cS^{\tor}_{X/k} \subseteq \cS^{\et}_{X/k}$. By Proposition~\ref{prop:toricbrauer}, it is enough to prove that, for every étale neighbourhood $(Y,r)$ and every finite extension $k'/k$, the relative Brauer group $\br(Y_{k'}/k')$ is trivial. Since $s$ is Selmer, $r$ is Selmer as well (the lift of a geometric section along a finite étale cover is again geometric by Proposition~\ref{prop:cartesian}); this implies that $Y(k_{\nu})$ is non-empty for every finite place $\nu$. Moreover, $Y(k_{\nu})$ is non-empty at infinite places as well, either because $k_{\nu}=\CC$ or by the real section conjecture \cite[Theorem 229, p.~208]{stix}. It follows that $\br(Y_{k'}/k')$ is trivial for every finite extension $k'/k$.

	For the converse, let $s\in\cS_{X/k}^{\tor} \subset \cS^{\et}_{X/k}$ be a toric Galois section. Fix a finite place $\nu$. We want to show that the image $s_{\nu}$ of $s$ in $\cS^{\et}_{X_{k_{\nu}}/k_{\nu}}$ is geometric. By Theorem~\ref{thm:torQp}, this is equivalent to showing that $s_{\nu} \in \cS^{\tor}_{X_{\nu}/k_{\nu}}$.
	
	Assume by contradiction that this is false. By Proposition~\ref{prop:torindex}, up to replacing $X$ with an étale neighbourhood we can assume that the index of $X_{k_{\nu}}$ is not $1$ (neighbourhoods of $s$, after base change, are cofinal among neighbourhoods of $s_{\nu}$, since their geometric base changes are cofinal among finite étale covers of $X_{\bar{k}}$), and hence there exists a Brauer class $b \in \br(X_{k_{\nu}}/k_{\nu})$ thanks to a theorem of Roquette and Lichtenbaum \cite{lichtenbaum}.
	
	Write $F = \bar{k}^{k_{\nu}}$. Since $\br(k) \to \br(k_{\nu})$ is surjective by the Brauer--Hasse--Noether theorem, then $\br(F) \to \br(k_{\nu})$ is surjective as well; let $P$ be a Brauer--Severi variety over $F$ such that $P_{k_{\nu}}$ is associated with $b$. By Lemma~\ref{lem:elemdes}, there exists a morphism $X_{F}\to P$, and hence $\br(X_{F}/F)$ is non-trivial.
	
	By Theorem~\ref{thm:basechange}, $\ptor_{X_{F}/F}=\ptor_{X/k}\times\spec F$, hence $s$ induces an element of $\cS_{X_{F}/F}^{\tor}$. This is in contradiction with Corollary~\ref{cor:brtor}.
\end{proof}

\section{The Tannakian interpretation of the toric fundamental group}\label{sect:cactus}

In this section we sketch a Tannakian interpretation of the toric fundamental gerbe. The rest of the paper is logically independent of this interpretation.

Throughout this section, $X$ is a smooth, projective, geometrically connected variety over a field $k$, of dimension $n\geq 1$. Fix an ample Cartier divisor $H$ on $X$. When $n=1$, the construction uses the ordinary degree of a line bundle and does not require a choice of $H$. We refer to \cite{deligne-milne} for the terminology of Tannakian categories.

\subsection{Pure cactuses and slopes}

If $f:Y\to X$ is a finite étale cover and $L$ is a line bundle on $Y$, write
\[d_H(L) = \deg_k \bigl(c_1(L)\cdot(f^*H)^{n-1}\bigr).\]
In particular, if $g:Z\to Y$ is a connected finite étale cover, projection formula gives
\[d_H(g^*L)=\deg(Z/Y) \cdot d_H(L).\]
Here and below the polarization on a cover is always pulled back from $X$.

\begin{definition}\label{def:pure-cactus}
	A vector bundle $V$ on $X$ is a \emph{pure cactus of slope $q\in\QQ$} if there exists a connected finite étale cover $f:Y\to X$ such that
	\[
	f^*V\simeq\bigoplus_{i=1}^{r}L_i,
	\qquad d_H(L_i)=\deg(Y/X) \cdot q\quad\text{for every }i.
	\]
	The cover and the splitting are not part of the datum.
	
	A \emph{cactus} on $X$ is a vector bundle with a specified finite $\QQ$-grading
	\[
	V=\bigoplus_{q\in\QQ}V_q,
	\]
	where $V_q$ is pure of slope $q$. A morphism of cactuses is a morphism of vector bundles preserving the grading. We denote this category by $\cact_{X,H}$, or simply $\cact_X$ when the polarization is understood.
\end{definition}

The slope $q$ of a pure cactus coincides with the slope of the underlying bundle: if $f: Y \to X$ is a finite étale cover with a splitting $f^*V\simeq\bigoplus_{i=1}^{r}L_i$, we get
\[q = \frac{d_{H}(L_{1})}{\deg(Y/X)} = \frac{(c_{1}(L_{1}) + \dots + c_{1}(L_{r}))\cdot(f^*H)^{n-1}}{r \cdot \deg(Y/X)} = \]
\[= \frac{c_{1}(f^{*}V) \cdot (f^*H)^{n-1}}{\rk(V) \cdot \deg(Y/X)} = \frac{c_{1}(V) \cdot H^{n-1}}{\rk(V)} = \mu_{H}(V).\]

The grading, however, is part of the datum of a general cactus, because we do not want a morphism of cactuses to mix different slopes. For example, on $\PP^{1}$ a non-zero section $\cO \to \cO(1)$ is not a morphism of cactuses, because $\cO$ has slope $0$ whereas $\cO(1)$ has slope $1$.

\begin{lemma}\label{lem:line-comp}
	Let $L,M$ be line bundles on a smooth, projective, connected $k$-scheme $Y$ of positive dimension, with equal degree with respect to some ample line bundle. Every nonzero morphism $L\to M$ is an isomorphism. 
\end{lemma}

\begin{proof}
	A nonzero morphism $L\to M$ gives a global section of $M\otimes L^{-1}$. Since $M \otimes L^{-1}$ has degree $0$ with respect to an ample line bundle, every non-zero section is never-vanishing. It follows that $L \to M$ is an isomorphism.
\end{proof}

\begin{proposition}\label{prop:cact-tannaka}
	The category $\cact_X$ is a $k$-linear Tannakian category. Its tensor product, dual and unit are given by
	\[
	(V\otimes W)_q=\bigoplus_{a+b=q}V_a\otimes W_b,
	\qquad
	(V^\vee)_q=(V_{-q})^\vee,
	\qquad
	\mathbbm{1}=\cO_X\text{ in slope }0.
	\]

	Every rational point $x\in X(k)$ defines a fibre functor $\omega_x(V)=V|_x$; in particular, $\cact_X$ is neutral if $X(k)\neq\emptyset$.
\end{proposition}

\begin{proof}
	Let us check that $\cact_{X}$ is abelian. Given a morphism of cactuses $V \to W$, we need to check that kernel and cokernel are cactuses. Clearly, we may assume that $V$, $W$ are pure of the same slope $q$. After passage to a finite étale cover, we might also assume that $V$ and $W$ are direct sums of line bundles of a fixed degree $d$. By Lemma~\ref{lem:line-comp}, we may furthermore reduce to the case $V = L^{\oplus a}$, $W = L^{\oplus b}$ for a fixed line bundle $L$. By Lemma~\ref{lem:line-comp} again, a morphism $L^{\oplus a} \to L^{\oplus b}$ is simply given by an $a \times b$ matrix. Kernel and cokernel of $V \to W$ are then isomorphic to $L^{\oplus a - r}$, $L^{\oplus b - r}$, where $r$ is the rank of the matrix. 

	The rest of the proof is straightforward and left to the reader.
\end{proof}

\subsection{The gerbe of cactuses}

For an object $V$ of a Tannakian category, write $\langle V\rangle$ for the smallest strictly full Tannakian subcategory containing $V$, and similarly $\langle V,W\rangle=\langle V\oplus W\rangle$. These finitely generated subcategories correspond to the finite type quotient gerbes of the gerbe of fibre functors.

\begin{lemma}\label{lem:cactus-pro-vt}
 The gerbe $\operatorname{Fib}(\cact_X)$ is pro-vt.
\end{lemma}

\begin{proof}
	Fix a cactus $V$, and choose a connected finite étale Galois cover $f:Y\to X$ with Galois group $G$ splitting its graded pieces into lines $L_i$. The bundle $A=f_*\cO_Y$ is finite, because $f^{*}A = f^{*}f_{*}\cO_{Y} \simeq \cO_{Y}^{|G|}$. In particular, it is a pure cactus of slope $0$.
	
	Choose a geometric point of $X$, let $\omega:\cact_{X} \to \vect_{\bar{k}}$ be the fibre functor and let $\aut(\omega)$ be the geometric Tannaka group of $\langle V,A\rangle$. The representation of $\aut(\omega)$ on $\omega(V\oplus A)$ is faithful, whereas the representation on $\omega(A)$ has finite image, because $A$ is a finite bundle. Let $H$ be the kernel of $\aut(\omega) \to \GL(\omega(A))$.

	We get a faithful action of $H$ on $\omega(V)$. Let $Y_{0}$ be the connected component of $Y_{\bar{k}}$ containing a chosen lift of the base point. By adjunction and the projection formula, pullback identifies $\langle V|_{Y_0}\rangle$ with the Tannakian quotient of $\langle V,A\rangle_{\bar{k}}$ by $\langle A\rangle_{\bar{k}}$ defined by the restriction of $\omega$. The Tannaka group of this quotient is $H$. Since $V|_{Y_0}$ is a direct sum of line bundles, $\omega(V)$ is a direct sum of one-dimensional $H$-representations, so $H$ acts diagonally and faithfully. This is sufficient to conclude that the connected component of $\aut(\omega)$ is of multiplicative type. Taking the limit along all $V$ proves that $\operatorname{Fib}(\cact_X)$ is pro-vt.
\end{proof}

The forgetful functor defines a canonical morphism $u : X \to \operatorname{Fib}(\cact_X)$. Since $\operatorname{Fib}(\cact_X)$ is pro-vt, the universal property of the toric fundamental gerbe gives a factorization $X \to \Pi^{\tor}_{X/k} \to \operatorname{Fib}(\cact_X)$.

\begin{theorem}\label{thm:cactus-gerbe}
	The morphism $\Pi^{\tor}_{X/k} \to \operatorname{Fib}(\cact_X)$ is an isomorphism.
\end{theorem}

\begin{proof}
	This can be checked after passing to the separable closure $k^{s}$ by Theorem~\ref{thm:basechange} and thanks to the fact that $\operatorname{Fib}(\cact_{X_{k^{s}}})$ identifies naturally with $\operatorname{Fib}(\cact_X) \otimes_{k} k^{s}$. Assume $k = k^{s}$.
	
	Let us show that $X \to \operatorname{Fib}(\cact_X)$ has the universal property of the toric fundamental gerbe. Since $k = k^{s}$ and $X$ is smooth, we may fix a $k$-point in $X$. Let $\Delta$ be a vt-gerbe with a morphism $h: X \to \Delta$, the base point identifies $\Delta = BG$ for some vt group $G$. We may replace $BG$ with the intermediate gerbe in the canonical factorization of $\ptor_{X/k}\to BG$ given by Proposition~\ref{prop:canonical} and assume that the induced morphism $\ptor_{X/k}\to BG$ is locally full. Write $T=G^0$. This corresponds to a connected $G$-torsor $P \to X$, a connected finite étale cover $f: Y = P/T \to X$ and a $T$-torsor $P \to Y$.
	
	We want to give a factorization $X \to \operatorname{Fib}(\cact_{X}) \to BG$; equivalently, we need to extend $\vect_{BG} \to \vect_{X}$ to an exact $k$-linear tensor functor $\operatorname{Rep}G \to \cact_{X}$.
	
	Write $\Omega$ for the group of characters of $T$. Since $T$ is abelian, we get a natural action of $G/T$ on $T$ by conjugation, and hence an action on $\Omega$. 
	
	For every character $\chi: T \to \GG_{m}$, we get a $\GG_{m}$-torsor induced by $P$ on $Y$, and hence a line bundle $L_{\chi}$ on $Y$. Define the slope 
	\[q_{\chi} = \frac{\deg_k\bigl(c_1(L_\chi)\cdot(f^*H)^{n-1}\bigr)}{\deg(Y/X)}.\]
	Notice that, if $\chi$ and $\chi'$ are in the same $G/T$ orbit, then $q_{\chi} = q_{\chi'}$ thanks to the fact that $H$ is defined over $X$.
	
	Let $W$ be a representation of $G$, corresponding to a vector bundle $\cW \to BG$. We have a splitting $W|_{T} = \bigoplus_{\chi \in \Omega} W_{\chi}$. For every orbit $S \subset \Omega$, write $W_{S} = \bigoplus_{\chi \in S} W_{\chi}$, it is immediate to check that $W = \bigoplus_{S}W_{S}$ as a $G$-representation. 
	
	Hence, we get a splitting $\cW = \bigoplus_{S} \cW_{S}$. Since
	\[f^*h^*\cW \simeq \bigoplus_\chi L_\chi\otimes_k W_\chi,\]
	we see that the pullback of $\cW_{S}$ is a pure cactus for every orbit $S$, hence $\cW$ has a natural structure of cactus. This gives the desired factorization $X \to \operatorname{Fib}(\cact_{X}) \to BG$, concluding the proof.
\end{proof}

Recall that a Tannakian category $\cC$ is algebraic if $\cC = \langle V \rangle$ for some $V \in \cC$, equivalently if its gerbe of fibre functors is algebraic \cite[Remark 3.10]{deligne-milne}.

\begin{definition}
 A \emph{quasi-fibre functor} of $\cC$ is a choice of an isomorphism class $[F_V]$ of $k$-valued fibre functors on $\langle V \rangle$ for every $V \in \cC$, such that
 \[
 [F_W|_{\langle V \rangle}]=[F_V]\quad\text{whenever }\langle V \rangle \subseteq \langle W \rangle.
 \]
\end{definition}

This is a compatibility condition on isomorphism classes; coherent choices of isomorphisms are not part of the datum. Restriction defines a natural map from isomorphism classes of fibre functors on $\cC$ to quasi-fibre functors. The set of quasi-fibre-functors identifies naturally with $\operatorname{Fib}(\cC)(k)_{\rmq}$. In particular, a toric Galois section in $\cS^{\tor}_{X/k} = \Pi^{\tor}_{X/k}(k)_{\rmq}$ corresponds to a quasi-fibre functor on $\cact_{X}$.

\printbibliography

\end{document}